\documentclass[arxiv]{agn_article}
\usepackage{agn_all}
\usepackage{float}
\usepackage{tabularx}
\usepackage{caption}
\graphicspath{{images/}}

\newcommand{\comsol}{\textit{Comsol Multiphysics}\textsuperscript{\textregistered}}
\newcommand{\mathematica}{\textit{Mathematica}}
\newcommand{\citet}[2][]{\citeauthor{#2} \cite[#1]{#2}}

\providecommand{\mulogiso}{\mu_{\rm log}^{\rm iso}}
\providecommand{\mueuclidiso}{\mu_{\rm euclid}^{\rm iso}}
\providecommand{\muiso}{\mu^{\rm iso}}
\providecommand{\kappaiso}{\kappa^{\rm iso}}
\providecommand{\kappalogiso}{\kappa_{\rm log}^{\rm iso}}
\providecommand{\kappaeuclidiso}{\kappa_{\rm euclid}^{\rm iso}}

\providecommand{\Ciso}{\C_{\rm iso}}
\providecommand{\Caniso}{\C_{\rm aniso}}
\providecommand{\mus}{\mu^*}

\defineauthor{panos}{Panos~Gourgiotis}{samps}{pgourgiotis@mail.ntua.gr}
\defineaffiliation{samps}{Mechanics Division, SAMPS, National Technical University of Athens, Zographou, GR-15773, Greece}
\defineauthor{lewintankit}{Peter~Lewintan}{agn,kit}{peter.lewintan@uni-due.de}
\defineaffiliation{kit}{Faculty of Mathematics, Karlsruhe Institute of Technology, Englerstraße 2, 76131 Karlsruhe, Germany}
\defineauthor{sky}{Adam~Sky}{lux}{adam.sky@uni.lu}
\defineaffiliation{lux}{Faculty of Science, Technology and Medicine, University of Luxembourg,  6 Av. de la Fonte, 4364 Esch-Belval Esch-sur-Alzette, Luxembourg}
\defineauthor{franz}{Franz~Gmeineder}{konstanz}{franz.gmeineder@uni-konstanz.de}
\defineaffiliation{konstanz}{Department of Mathematics and Statistics, University of Konstanz, Universitätsstrasse 10, 78457 Konstanz, Germany}

\begin{document}

\title{Yet another best approximation isotropic elasticity tensor\\ in plane strain}
\knownauthors[vossd]{vossd,panos,lewintankit,sky,neff}
\date{\today}
\maketitle

\begin{abstract}
\noindent
For plane strain linear elasticity, given any anisotropic elasticity tensor $\C_{\rm aniso}$, we determine a best approximating isotropic counterpart $\Ciso$.
This is not done by using a distance measure on the space of positive definite elasticity tensors (Euclidean or logarithmic distance) but by considering two simple isotropic analytic solutions (center of dilatation and concentrated couple) and best fitting these radial solutions to the numerical anisotropic solution based on $\C_{\rm aniso}$.
The numerical solution is done via a finite element calculation, and the fitting via a subsequent quadratic error minimization.
Thus, we obtain the two Lam\'e-moduli $\mu$, $\lambda$ (or $\mu$ and the bulk-modulus $\kappa$) of $\Ciso$.
We observe that our so-determined isotropic tensor $\Ciso$ coincides with neither the best logarithmic fit of Norris nor the best Euclidean fit. Our result calls into question the very notion of a best-fit isotropic elasticity tensor to a given anisotropic material.
\end{abstract}

\\[.65em]
\noindent\textbf{AMS 2010 subject classification:
	74A10, 
	74B05, 
	74M25  
}
{\parskip=-0.4mm\tableofcontents}

%
%
%
\section{Introduction}
Linear elasticity theory is a mature science that allows us to predict with sufficient accuracy the response of materials or structures under incremental loading once their constitutive tensors are known. In practice, however, a problem resides with determining the elasticity tensor from measurements, and sometimes the symmetry classes are too involved and call for some simplification.
The simplification might consider an elasticity tensor with more symmetries, replacing one with less symmetries. It is clear that an elasticity tensor with larger symmetries has fewer free parameters to be determined. A maximum simplification occurs if one wants to replace a given symmetry class with its isotropic counterpart - if such an object is well-defined at all.

There are several procedures to be found in the literature aspiring to find this isotropic counterpart. It is well known that the "isotropic" result depends heavily on the employed objective function that needs to be minimized to obtain the "best-fit" isotropic tensor. Proposed are, e.g. taking the Euclidean distance on the space of elasticity tensors \cite{kochetov2009estimating,cavallini1999best}, the logarithmic distance \cite{norris2006isotropic,moakher2006closest}, or using the Euler--Lagrange method \cite{antonelli2022distance,azzi2023distance}. Alternatively, the anisotropic acoustic tensor $\Caniso.q\otimes q$ may be interpreted as a $q$-direction dependent quadratic form and
approximated in some sense, e.g. via a "slowness fit" \cite{kochetov2009estimating}: find the best fit in terms of circles (representatives of the isotropic elasticity tensor $\Ciso$) to a non-circular wavefront (generated via $\Caniso$) \cite{danek2016closest,danek2018effects}. In all these cases, a direct measure acting somehow on $\Caniso$ is considered. 

Even for the Euclidean distance, surprisingly, the results differ if a representation of the elasticity tensor in $6\times6$ Voigt notation is used or the full $9\times9$ component representation is taken \cite{danek2016closest,danek2018effects}. Another observation is that it should not matter whether we use the elasticity tensor $\Caniso$, or the corresponding compliance tensor $\C_{\rm aniso}\inv$ for the determination of the best isotropic fit $\Ciso$. This requirement already excludes the simple approach by the Euclidean distance. On the other hand, the geodesic approach \cite{agn_martin2014minimal} indeed has the significant advantage of delivering the same results, i.e. the geodesic distances satisfy the basic invariance property
\begin{equation}
	\dist_{\rm geod}(\C_{\rm iso},\C_{\rm aniso})=\dist_{\rm geod}(\C_{\rm iso}\inv,\C_{\rm aniso}\inv).
\end{equation}

Here, we propose yet another method to obtain a best-fit isotropic tensor, albeit restricted to plane strain. We consider two representative isotropic analytical solutions in plane strain elasticity (the concentrated couple and the center of dilatation). Geometrically, these solutions exhibit straightforward radial symmetry due to planar isotropy. Both solutions can be also obtained numerically by applying the corresponding traction on a small but finite circle around the origin.  To connect this to the anisotropic setting, we take a given planar $\Caniso$ (mostly cubic) and subject it to the same traction on the same inner circle. The finite element solution found computationally, here done with \comsol, is no longer radial. We now best fit this non-radial solution with the analytical radial solution and are thus able to determine $\Ciso$ uniquely.

We investigate two scenarios: The first considers only the norm of the displacement vector $\norm{u(x_1,x_2)}$ for the ensuing best fit. The second considers the full displacement field $u(x_1,x_2)$ instead. Both methods are compared in Section \ref{sec:discussion}. They may greatly differ from either the Euclidean or logarithmic approach depending on the strength of the anisotropy, as indeed all discussed approaches deliver the correct isotropic result for an isotropic input. We believe that for our academic plane strain problem the full-field resolution and optimization represent a valuable alternative to more direct calculations.  The computational burden is manageable because the FEM calculation has to be done only twice, once for the center of dilatation and once for the concentrated couple, and is applicable to any symmetry class.

The paper is structured as follows. We derive the linear elastic solution for a concentrated couple and center of dilatation for planar cubic materials using Green's functions and recall the well-known isotropic linear elastic case. We show that they can straightforwardly be used to correctly retrieve isotropic parameters via a full-field simulation in \comsol\, and quadratic error minimization (here done in \mathematica) for both the norm-scenario and the full displacement approach. Further, our results reveal the necessary domain size for sufficient accuracy. Since both solutions decay to zero far away from the center, it becomes inconsequential whether one tries to find a matching amplitude for an additional Dirichlet boundary displacement, or simply fixes the displacement field to zero on the entire boundary.

Next, we apply the same setting to a given cubic tensor $\Caniso$. The finite element solution for a center of dilatation and a concentrated couple is obtained, and subsequently best-fitted against the analytical isotropic solution, thus determining $\Ciso$. We repeat this procedure for several distinct cubic tensors and compare the obtained values with the ones stemming from calculations using Euclidean and logarithmic minimization. Finally, we discuss conclusions and outlook.
%
%
%
\subsection{Best approximated isotropic tensor by Norris - logarithmic distance}\label{sec:norris}
The general problem of finding the closest isotropic elasticity tensor $\C_{\rm iso}$ for a given anisotropic elasticity tensor $\C_{\rm aniso}$ has intrigued the scientific community for some time. Here, we follow the approach and complete solution provided by A.\ Norris \cite{norris2006isotropic}. The first observation is that it should not matter whether we use the elasticity tensor $\C$ or its inverse, the compliance tensor $\C\inv$. This requirement already excludes the simple approach by the Euclidean distance. Instead one has to revert to the logarithmic distance which satisfies the basic invariance property
\begin{equation}
	\dist_{\rm log}(\C_{\rm iso},\C_{\rm aniso})=\dist_{\rm log}(\C_{\rm iso}\inv,\C_{\rm aniso}\inv).
\end{equation}

We give a summary of the best approximated isotropic tensors by Norris \cite{norris2006isotropic} for both the Euclidean and the logarithmic distance. For this, we show a short calculation leading to the explicit formula that allows us to compute the corresponding closest isotropic tensor $\Ciso$ for any anisotropic elasticity tensor $\Caniso$ with cubic symmetry. This specific type of anisotropic tensor appears frequently
especially in the setting of mechanical metamaterials \cite{agn_madeo2015wave,agn_agostino2019dynamic,agn_voss2023modeling}.

Norris identification varies with the norm applied (Euclidean or logarithmic)
\begin{align}
	\dist_{\rm euclid}(A,B)=\norm{A-B}\,,\qquad\dist_{\rm log}(A,B)=\norm{\log A-\log B}.
\end{align}
In the context of elasticity, using Voigt notation, we deal with positive definite matrices $\C\in\Symp(6)$. Here, we consider materials with isotropic and cubic symmetry
\begin{align}
	\C_{\rm iso}=\matr{
		\kappa+\mu & \kappa-\mu & \kappa-\mu & 0 & 0 & 0\\
		\kappa-\mu & \kappa+\mu & \kappa-\mu & 0 & 0 & 0\\
		\kappa-\mu & \kappa-\mu & \kappa+\mu & 0 & 0 & 0\\
		0 & 0 & 0 & \mu & 0 & 0\\
		0 & 0 & 0 & 0 & \mu & 0\\
		0 & 0 & 0 & 0 & 0 & \mu
	},\qquad\C_{\rm cubic}=\matr{
		\kappa+\mu & \kappa-\mu & \kappa-\mu & 0 & 0 & 0\\
		\kappa-\mu & \kappa+\mu & \kappa-\mu & 0 & 0 & 0\\
		\kappa-\mu & \kappa-\mu & \kappa+\mu & 0 & 0 & 0\\
		0 & 0 & 0 & \mus & 0 & 0\\
		0 & 0 & 0 & 0 & \mus & 0\\
		0 & 0 & 0 & 0 & 0 & \mus
	}.\label{eq:voigtNotation}
\end{align}
Thus, in the isotropic case (for plane strain), we have two independent material parameters: the bulk modulus $\kappa>0$ and the shear modulus $\mu>0$. For the cubic case we also have a third independent material parameter: the second shear modulus $\mus>0$.

We recall the well-known Euclidean/logarithmic distance between two cubic matrices $\C_{\rm cubic}(\kappa_1,\mu_1,\mu_1^*)$ and  $\C_{\rm cubic}(\kappa_2,\mu_2,\mu_2^*)$, cf. \cite[eq. (22)]{norris2006isotropic}
\begin{align}
	\dist_{\rm euclid}(\C_{\rm cubic}(\kappa_1,\mu_1,\mu_1^*),\C_{\rm cubic}(\kappa_2,\mu_2,\mu_2^*))&=\sqrt{9\.(\kappa_1-\kappa_2)^2+8\.(\mu_1-\mu_2)^2+12\.(\mu_1^*-\mu_2^*)^2},\\
	\dist_{\rm log}(\C_{\rm cubic}(\kappa_1,\mu_1,\mu_1^*),\C_{\rm cubic}(\kappa_2,\mu_2,\mu_2^*))&=\sqrt{\log^2\frac{\kappa_1}{\kappa_2}+2\.\log^2\frac{\mu_1}{\mu_2}+3\.\log^2\frac{\mu_1^*}{\mu_2^*}}.
\end{align}
Looking at \eqref{eq:voigtNotation}, we can write every isotropic matrix $\C_{\rm iso}(\kappa,\mu)$ as the corresponding cubic matrix $\C_{\rm cubic}(\kappa,\mu,\mu)$.
Next, we want to find the closest isotropic matrix $\C_{\rm iso}(\kappaiso,\muiso)$ for any cubic matrix $\C_{\rm cubic}(\kappa,\mu,\mu^*)$. For the Euclidean case, we compute
\begin{align}
	\min_{\kappaiso,\muiso}\.\dist_{\rm euclid}(\C_{\rm iso}(\kappaiso,\muiso),\C_{\rm cubic}(\kappa,\mu,\mu^*))&\overset{\phantom{\kappaiso=\kappa}}{=}\min_{\kappaiso,\muiso}\.\sqrt{9\.(\kappaiso-\kappa)^2+8\.(\muiso-\mu)^2+12\.(\muiso-\mu^*)^2}\notag\\
	&\overset{\kappaiso=\kappa}{=}\min_{\muiso}\.2\.\sqrt{2\.(\muiso-\mu)^2+3\.(\muiso-\mu^*)^2},
\end{align}
such that $\kappaiso=\kappa$. For the optimal $\muiso$, we compute the roots of the first derivative
\begin{align}
    0&=\dd{\muiso}\left[2\.(\muiso-\mu)^2+3\.(\muiso-\mu^*)^2\right]=10\.\muiso-4\mu-6\.\mu^*\qquad\iff\qquad\mueuclidiso=\frac25\.\mu+\frac35\.\mu^*.\label{eq:norrisEuclid}
\end{align}
Thus, for any given cubic matrix $\C_{\rm cubic}(\kappa,\mu,\mu^*)$, its best approximated isotropic matrix using the Euclidean distance is $\C_{\rm iso}(\kappa,\frac25\.\mu+\frac35\.\mu^*)$. 

We repeat the procedure with the logarithmic distance as well
\begin{align}
	\min_{\kappaiso,\muiso}\.\dist_{\rm euclid}(\C_{\rm iso}(\kappaiso,\muiso),\C_{\rm log}(\kappa,\mu,\mu^*))&\overset{\phantom{\kappaiso=\kappa}}{=}\min_{\kappaiso,\muiso}\.\sqrt{\log^2\frac{\kappa}{\kappaiso}+2\.\log^2\frac{\muiso}{\mu}+3\.\log^2\frac{\muiso}{\mu^*}}\notag\\
	&\overset{\kappaiso=\kappa}{=}\min_{\muiso}\.\sqrt{2\.\log^2\frac{\muiso}{\mu}+3\.\log^2\frac{\muiso}{\mu^*}}.
\end{align}
Again $\kappaiso=\kappa$, and for the optimal $\muiso$, we compute the roots of the first derivative
\begin{align}
	&&0&=4\.\log\frac{\muiso}{\mu}\.\frac{1}{\muiso}+6\.\log\frac{\muiso}{\mu^*}\.\frac{1}{\muiso}=\frac2\muiso\left(5\.\log\muiso-2\.\log\mu-3\.\log\mu^*\right)\notag\\
	&\iff&\log\bigl((\muiso)^5\bigr)&=\log\bigl(\mu^2\bigr)+\log\bigl((\mu^*)^3\bigr)\qquad\iff\mulogiso=\sqrt[5]{\mu^2(\mu^*)^3}.\label{eq:norrisLog}
\end{align}
As such, for any given cubic matrix $\C_{\rm cubic}(\kappa,\mu,\mu^*)$, its best approximated isotropic matrix using the logarithmic distance is $\C_{\rm iso}(\kappa,\sqrt[5]{\mu^2(\mu^*)^3})$. 
%
%
%
%
\subsubsection{Reverse formula}
The best approximation formulas \eqref{eq:norrisEuclid} and \eqref{eq:norrisLog} can also be used to find an inverted relation, i.e. indicate which type of cubic matrices have the same best approximated isotropic matrix.
In the formulation of Norris, the bulk modulus $\kappa=\kappaiso$ remains always constant regardless of the employed distance measure.
However, different pairs of $(\mu,\mu^*)$ can lead to identical values $\mueuclidiso$, $\mulogiso$ when either the Euclidean distance \eqref{eq:norrisEuclid} or the logarithmic distance \eqref{eq:norrisLog} is used.

For the Euclidean distance we have the linear relation $\mueuclidiso=\frac25\.\mu+\frac35\.\mu^*$ and thus we can make the ansatz $\C_{\rm cubic} 
(\kappa,\mueuclidiso+3\.c,\mueuclidiso-2\.c)$ with the free parameter $c\in\R$. We verify
\begin{equation}
	\frac25\.(\mueuclidiso+3\.c)+\frac35\.(\mueuclidiso-2\.c)=\muiso+c-c=\mueuclidiso\,.
\end{equation}
For the logarithmic distance $\mulogiso=\sqrt[5]{\mu^2(\mu^*)^3}$, we choose the ansatz $\C_{\rm cubic}\bigl(\kappa,\mulogiso\.c^3,\frac{\mulogiso}{c^2}\bigr)$ with the free parameter $c\in\Rp$. We affirm
\begin{equation}
	\sqrt[5]{(\mulogiso\.c^3)^2\left(\frac{\mulogiso}{c^2}\right)^3}=\sqrt[5]{(\mulogiso)^5\frac{c^6}{c^6}}=\mulogiso\,.
\end{equation}
\begin{figure}[h!]
	\centering
	\includegraphics[width=0.5\linewidth]{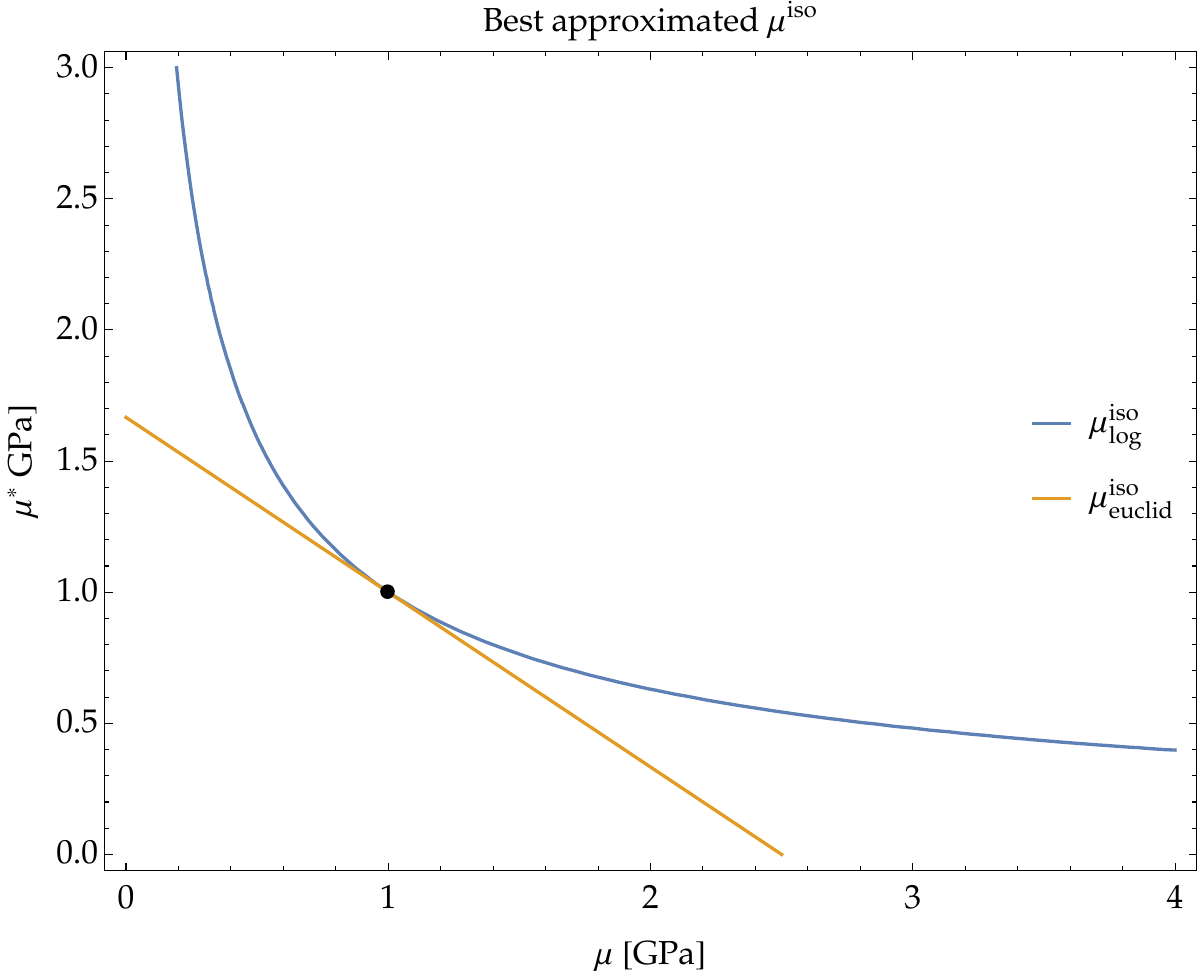}\centering
	\caption{Visualization of Norris formula for the Euclidean distance (orange line) and the logarithmic distance (blue line). The contour plot shows the associated best approximated $\muiso=1$ for given $\mu,\mu^*$ (cubic symmetry) while the bulk modulus $\kappa=\kappaiso$ has no effect. Both approximations coincide for $\mu=\mu^*=1$.}
\end{figure}
%
%
%
\section{Green's functions for planar cubic materials}\label{sec:GreenFunction}
In this section, we derive in closed form the plane strain Green's functions for a cubic material. Employing the obtained Green's functions, we construct the solutions, first for the concentrated couple, and then for the center of dilatation. The solutions are subsequently used in our best approximation procedure of the isotropic elasticity tensor. 
The analysis that will be followed is based on double Fourier transforms and is analogous to the one presented in Bigoni and Gourgiotis \cite{BG16} for orthotropic couple stress materials. 

For a cubic material under plane strain conditions, the displacement equations of equilibrium read
\begin{equation}
\label{EQs}
\begin{aligned}
\left( \lambda +2\mu  \right)u_{1,x_1 x_1}+\mu^*\.u_{1,x_2x_2}+\left( \lambda +\mu^*\right)u_{2,x_1x_2}+{{f}_{1}}=0\,,\\
\left( \lambda +2\mu  \right)u_{2,x_2 x_2}+\mu^*\.u_{2,x_1 x_1}+\left( \lambda +\mu^* \right)u_{1,x_1x_2}+{{f}_{2}}=0\,,
\end{aligned}
\end{equation}
where $u_i \equiv u_i(x_1,x_2)$ are the displacement components and $f_i \equiv f_i(x_1,x_2)$ are the components of the body force. The strain energy density is positive definite in the plane strain case examined here when the material moduli satisfy the following inequalities
\begin{equation}
\mu>0, \qquad \lambda+\mu>0, \qquad \mu^*>0\,.
\end{equation}
The infinite body plane strain Green’s function for a concentrated force $f_i=P_i\,\delta(x_1)\delta(x_2)$ is obtained by employing the double exponential Fourier transform. The direct and inverse double Fourier transforms are defined as
\begin{align}
	\label{FT}
	&\hhat(\xi_1, \xi_2)=\int_{-\infty }^{+\infty }\int_{-\infty }^{+\infty} h(x_1,x_2)\. e^{i\.\xi\cdot x}\,\dx_1\dx_2\,, \\
	& h(x_1,x_2)=\frac{1}{4 \pi ^2}\int_{-\infty }^{+\infty }\int_{-\infty }^{+\infty} \hhat(\xi_1,\xi_2)\. e^{-i\.\xi\cdot x}\, \intd\xi_1 \intd\xi_2\,,
\end{align}
where $x=(x_1,x_2)$ and $\xi=(\xi_1,\xi_2)$ is the Fourier vector.
Applying the direct transform \eqref{FT} to the field equations \eqref{EQs} yields the following solution for the displacement field 
\begin{equation}
\label{Sol}
u_q(x_1,x_2)=\frac{P_p}{4 \pi^2} \int_{-\infty }^{\infty}\int_{-\infty }^{\infty} \frac{C_{pq}(\xi)}{D (\xi)} \,e^{-i\,\xi\cdot x}\, \intd\xi_1 \intd\xi_2 \, , \qquad q=1,2\,,
\end{equation}
where the components of $C_{pq}(\xi)$ are defined as
\begin{equation}
\begin{aligned}
&C_{11}(\xi)=\mu^*\xi_1^2+(\lambda +2\mu)\.\xi_2^2 \, ,\\
&C_{12}(\xi)=C_{21}(\xi)=-(\lambda+\mu^*)\. \xi_1 \xi_2 \, ,\\
&C_{22(\xi)}=(\lambda +2\mu)\.\xi_1^2+\mu^*\.\xi_2^2
\end{aligned}
\end{equation}
and 
\begin{equation}
D\left(  \xi  \right)=C_{11}C_{22}-C_{12}^2=\left( \lambda +2\mu  \right)\mu^*\.\xi _{1}^{4}+2\left( 2\mu \left( \lambda +\mu  \right)-\lambda \mu^* \right)\xi_1^2\xi_2^2+\left( \lambda +2\mu  \right)\mu^*\.\xi_2^4 \, .
\end{equation}
Note that $C_{pq}(\xi)$ is actually the cofactor of the acoustic tensor and $D(\xi)$  is the determinant of the acoustic tensor for a cubic plane strain material.

For a fixed value of the transformed variable $\xi_1\in\R$, the determinant $D(\xi)$ is a quartic non-homogeneous polynomial of the variable $\xi_2$ which has no real roots when the material is positive definite. In fact, the four roots of the characteristic quartic polynomial can be either purely imaginary or complex conjugate which are placed symmetrically in the four quadrants of the $\xi_2$-plane. The characteristic polynomial can be now written as
\begin{equation}
D( \xi)=\mu^*(\lambda +2\mu) \prod_{m=1}^2 \left( \xi_2-\xi_2^{(m)} \right) \left( \xi_2-\overline{\xi}_2^{\,(m)} \right),\qquad m=1,2\,, 
\end{equation}
with $\xi_2^{(m)}=i s_m \abs{\xi_1\!\.}$ and $\overline{\xi}_2^{(m)}$ being the complex conjugates of $\xi_2^{(m)}$. Also, 
\begin{equation}
s_{1,2}=\frac{1}{\sqrt{2 \mu^*(\lambda+2\mu)}} \sqrt{-p \mp \sqrt q}\,,
\end{equation}
where
\begin{equation}
p=2\lambda \mu^*-4\mu(\lambda +\mu) \, , \qquad q=16\mu \left( \mu -\mu^* \right)\left( \lambda +\mu  \right)\left( \lambda +\mu +\mu^* \right).
\end{equation}
It should be noted that for a fixed value of $\xi_1$ we have that $\text{Im}[\xi_2^{(m)}]>0$, while $\text{Im}[\bar{\xi}_2^{\,(m)}]<0$. Applying the residual theorem in conjunction with Jordan's
lemma, the integration with respect to $\xi_2$ in \eqref{Sol} yields a summation of residues of
poles at $\xi_2=\xi_2^{(m)}$ when $x_2<0$, and at $\xi_2=\bar{\xi}_2^{\,(m)}$ when $x_2>0$. In particular, for $x_2>0$, the original integration path running along the real axis in the $\xi_2$-plane is replaced by a closed contour taken in the lower $\xi_2$-plane so that the integrand decays as
$\lvert \xi_2 \rvert \to \infty$. In view of the above, the following relation is obtained 
\begin{equation}
\int_{-\infty }^{\infty} \frac{C_{pq}(\xi)}{D (\xi)} \,e^{-i\,\xi_2 x_2} \, \intd\xi_2=Q_{pq}(\xi_1,x_2) \, ,
\end{equation}
with
\begin{equation}
Q_{pq}(\xi_1,x_2)=-2 \pi i \sum_{m=1}^{2} \left[ \frac{C_{pq}(\xi) \, e^{-i\. \xi_2 x_2}}{\partial_{\xi_2} D(\xi)} \right]_{\xi_2=\bar{\xi}_2^{\,(m)}}.
\end{equation}
More specifically, we have that
\begin{equation}
\begin{aligned}
&Q_{11}=\frac{\pi}{\sqrt q \abs{\xi_1\!\.}}\,\left(Y_2^{(1)}  e^{-s_2 \abs{\xi_1\!\.} x_2}-Y_1^{(1)} e^{-s_1 \abs{\xi_1\!\.} x_2} \right),  \\
&Q_{12}=Q_{21}=i\frac{\pi}{q^{1/2} \, \xi_1}\, \left(\lambda+\mu^*\right) \left(e^{-s_1 \abs{\xi_1\!\.} x_2} -e^{-s_2 \abs{\xi_1} x_2\!\.}\right), \\
&Q_{22}=\frac{\pi}{\sqrt q\abs{\xi_1}}\,\left(Y_2^{(2)}  e^{-s_2 \abs{\xi_1\!\.} x_2}-Y_1^{(2)} e^{-s_1 \abs{\xi_1\!\.} x_2} \right), 
\end{aligned}
\end{equation}
with
\begin{equation}
Y_m^{(1)}=\frac{(\lambda+2\mu)s_m^2-\mu^*}{s_m} \, , \qquad Y_m^{(2)}=\frac{\mu^* s_m^2-(\lambda+2\mu)}{s_m} \, , \qquad m=1,2\,.
\end{equation}
Next, we utilize the fact that the components of $Q_{pq}$  are even functions of $\xi_1$ when $p=q$, and odd functions when  $p \ne q$, and employ the well-known results for the Fourier integrals
\begin{align}
\int_0^{\infty} \frac{e^{-s_m \xi_1 x_2}}{\xi_1} \, \cos{(\xi_1 x_1)} \, \intd\xi_1=\gamma -\log\sqrt{x_1^2+s_m^2 x_2^2}\, ,
\\
\int_0^{\infty} \frac{e^{-s_m \xi_1 x_2}}{\xi_1} \, \sin{(\xi_1 x_1)} \, \intd\xi_1=\frac{\pi}{2}-\tan^{-1}\left(\frac{s_m x_2}{x_1}\right),\notag
\end{align}
where the first integral is to be interpreted in the finite part sense and $\gamma\approx0.577$ is Euler's constant.
Neglecting terms corresponding to rigid body displacements in the $x_1x_2-$plane, we derive the final expressions for the  displacements due to a concentrated force:
\begin{itemize}
    \item \textit{horizontal force $P_1$}
    \begin{equation}
        u_1^{(P_1)}=\frac{P_1}{2\pi \.\sqrt q} \left(Y_1^{(1)}\log r_1-Y_2^{(1)}\log r_2\right), \qquad  u_2^{(P_1)}=\frac{P_1 (\lambda+\mu^*)}{2\pi \. \sqrt q} \left(\theta_2-\theta_1 \right),
    \end{equation}
    \item \textit{vertical force $P_2$}
    \begin{equation}
        u_1^{(P_2)}=\frac{P_2 (\lambda+\mu^*)}{2\pi \.\sqrt q} \left(\theta_2-\theta_1 \right), \qquad  u_2^{(P_2)}=\frac{P_2}{2\pi \.\sqrt q} \left(Y_1^{(2)}\log r_1-Y_2^{(2)}\log r_2\right),
    \end{equation}
    where
    \begin{equation}
        r_m=\sqrt{x_1^2+s_m^2 x_2^2}, \qquad \theta_m=\tan^{-1}\left(\frac{s_m x_2}{x_1} \right), \qquad m=1,2\,.
    \end{equation}
\end{itemize}
\begin{figure}[h!]
	\pgfmathsetmacro{\r}{1.5}
	\begin{minipage}{0.5\textwidth}
		\centering
		\begin{tikzpicture}
			\draw [->] (-1.7*\r,0) -- (1.7*\r, 0) node[below] {$x_1$};
			\draw [->] (0,-1.7*\r) -- (0, 1.7*\r) node[right] {$x_2$};
			\draw [->,red,ultra thick] (0.5*\r,-0.5*\r) -- (0.5*\r,0.5*\r) node[right,pos=0.8] {1};
			\draw [->,red,ultra thick] (0.5*\r,0.5*\r) -- (-0.5*\r,0.5*\r) node[above,pos=0.8] {1};
			\draw [->,red,ultra thick] (-0.5*\r,0.5*\r) -- (-0.5*\r,-0.5*\r) node[left,pos=0.8] {1};
			\draw [->,red,ultra thick] (-0.5*\r,-0.5*\r) -- (0.5*\r,-0.5*\r) node[below,pos=0.8] {1};
			\draw [<->,thick] (0,-0.1*\r) -- (0.5*\r,-0.1*\r) node[above,pos=0.5] {$d/4$};
			\draw [<->,thick] (-0.1*\r,0) -- (-0.1*\r,0.5*\r);
		\end{tikzpicture}	
	\end{minipage}
	\begin{minipage}{0.5\textwidth}
		\centering
		\begin{tikzpicture}
			\draw [->] (-1.7*\r,0) -- (1.7*\r, 0) node[below] {$x_1$};
			\draw [->] (0,-1.7*\r) -- (0, 1.7*\r) node[right] {$x_2$};
			\draw [->,red,ultra thick] (0.5*\r,0) -- (1.5*\r,0) node[below,pos=0.8] {1};
			\draw [->,red,ultra thick] (0,0.5*\r) -- (0,1.5*\r) node[right,pos=0.8] {1};
			\draw [->,red,ultra thick] (-0.5*\r,0) -- (-1.5*\r,0) node[below,pos=0.8] {1};
			\draw [->,red,ultra thick] (0,-0.5*\r) -- (0,-1.5*\r) node[right,pos=0.8] {1};
			\draw [<->,thick] (0,-0.1*\r) -- (0.5*\r,-0.1*\r) node[above,pos=0.5] {$d/4$};
			\draw [<->,thick] (-0.1*\r,0) -- (-0.1*\r,0.5*\r);
		\end{tikzpicture}
	\end{minipage}
	\caption{(Left) Schematic representation of the concentrated couple. (Right) Schematic representation of the center of dilatation.}\label{fig:DM}
\end{figure}
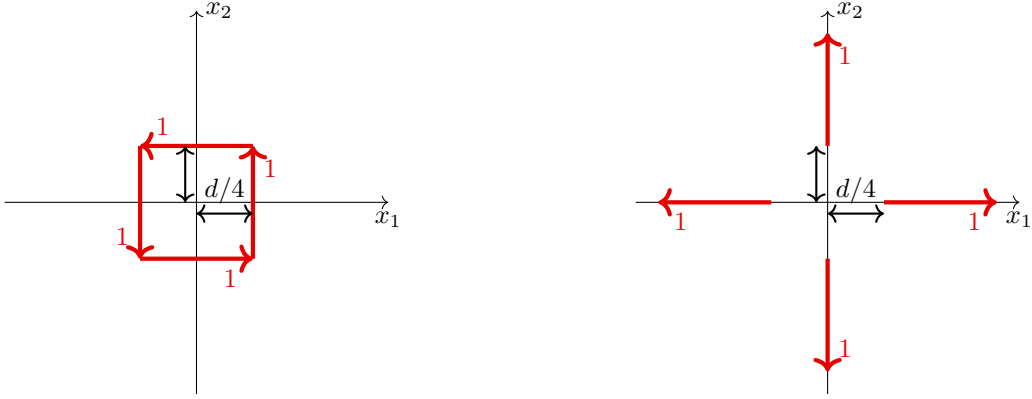
%
%
%
\subsection{Concentrated couple}
We will now use Green's function for a concentrated force to construct the solution for the concentrated couple in a cubic material.
Therein, the concentrated couple solution is constructed by superimposing two unit force dipoles with the moment as shown in Figure \ref{fig:DM}b). In this case, the following displacement field is derived
\begin{align}
    \label{M1}
    &u_1^{(M)}(x_1,x_2)=u_1^{(P_1)}(x_1,x_2+\frac{d}{4})-u_1^{(P_1)}(x_1,x_2-\frac{d}{4})-u_1^{(P_2)}(x_1+\frac{d}{4},x_2)+u_1^{(P_2)}(x_1-\frac{d}{4},x_2),\\
    &u_2^{(M)}(x_1,x_2)=u_2^{(P_1)}(x_1,x_2+\frac{d}{4})-u_2^{(P_1)}(x_1,x_2-\frac{d}{4})-u_2^{(P_2)}(x_1+\frac{d}{4},x_2)+u_2^{(P_2)}(x_1-\frac{d}{4},x_2).\notag
\end{align}
The solution for the concentrated couple of unit strength is then derived by taking the series expansion of \eqref{M1} with respect to the dipole arm $d$ and retaining only the linear term, which is equivalent to $\lim\limits_{d \to 0}\dd d u^{(D)}_j$, ($j=1,2$). In light of the above, the displacement field for a cubic material assumes the final form 
\begin{align}
    \label{eq:CCsolution}
    &u_1=-\frac{x_2}{4 \pi q^{1/2}} \, \left((\lambda+\mu^*)\left(\frac{s_1}{r_1^2}-\frac{s_2}{r_2^2}\right)-\frac{\left(Y_1^{(1)}-Y_2^{(1)}\right)s_1^2}{r_1^2} \right)-\frac{Y_2^{(1)}}{4\pi \mu^* (\lambda+2\mu)}\frac{x_1^2 x_2}{ r_1^2 r_2^2}\,,\\
    &u_2=\frac{x_1}{4 \pi q^{1/2}} \, \left(\frac{(\lambda+\mu^*)s_2+Y_2^{(2)}}{r_2^2}-\frac{(\lambda+\mu^*)s_1+Y_1^{(2)}}{r_1^2} \right).\notag
\end{align}
It is worth noting that a concurrent rigid rotation of the dipoles in Figure \ref{fig:DM} (Left) does not affect the above solution. 
We emphasize that the analytical solution \eqref{eq:CCsolution} for the concentrated couple for planar cubic materials is entirely algebraic, depending only on the three material parameters $\mu$, $\mu^*$ and $\lambda$ through $s_1$, $s_2$, $Y_1^{(1)}$ and $Y_2^{(1)}$.
In the isotropic case ($\mu^*=\mu$), the solution for the isotropic concentrated couple simplifies to
\begin{equation}
	u_1=-\frac{x_2}{4 \pi \mu\. r^2} \, , \qquad
	u_2=\frac{x_1}{4 \pi \mu\. r^2}\, .\label{eq:CClinear}
\end{equation}
\begin{figure}[h!]
	\centering
	\includegraphics[scale=0.7]{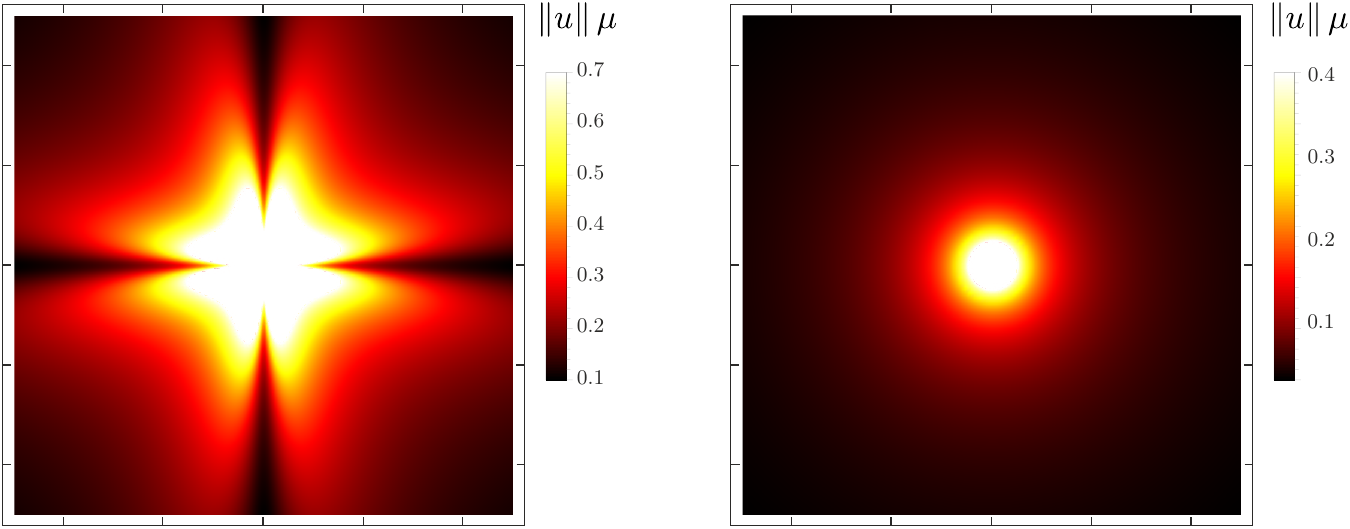}
	\caption{Variation of the normalized modulus of displacement $\norm{u} \mu$ for a concentrated couple (unit strength), (Left) cubic material ($\mu^*=0.1\. \mu$, $\lambda=0.29\. \mu$), (Right) isotropic material ($\lambda=0.29\. \mu$).}
	\label{dis}
\end{figure}
%
%
%
\subsection{Center of dilatation}
We now use Green's function for a concentrated force to construct the solution for the center of dilatation in a cubic material.
To this end, we superimpose two infinitesimal unit force dipoles at right angles that do not produce any moment as shown in Figure \ref{fig:DM} (Right), obtaining 
\begin{align}
\label{D1}
&u_1^{(D)}(x_1,x_2)=u_1^{(P_1)}(x_1-\frac{d}{4},x_2)-u_1^{(P_1)}(x_1+\frac{d}{4},x_2)+u_1^{(P_2)}(x_1,x_2-\frac{d}{4})-u_1^{(P_2)}(x_1,x_2+\frac{d}{4})\,,\\
&u_2^{(D)}(x_1,x_2)=u_2^{(P_1)}(x_1-\frac{d}{4},x_2)-u_2^{(P_1)}(x_1+\frac{d}{4},x_2)+u_2^{(P_2)}(x_1,x_2-\frac{d}{4})-u_2^{(P_2)}(x_1,x_2+\frac{d}{4})\,.\notag
\end{align}
The solution for the center of dilation is now
derived by taking the series expansion of \eqref{D1} with respect to the dipole arm $d$, 
retaining only the linear term which is equivalent to $\lim\limits_{d \to 0}\dd d u^{(D)}_j$, ($j=1,2$). In light of the above, the displacement field for a cubic material assumes the final form 
\begin{align}
\label{eq:CDsolution}
&u_1=\frac{x_1}{4 \pi\.\sqrt q} \, \left(\frac{(\lambda+\mu^*)s_1-Y_1^{(1)}}{r_1^2}-\frac{(\lambda+\mu^*)s_2-Y_2^{(1)}}{r_2^2} \right),
\\
&u_2=\frac{x_2}{4 \pi\.\sqrt q} \, \left((\lambda+\mu^*)\left(\frac{s_2}{r_2^2}-\frac{s_1}{r_1^2}\right)-\frac{\left(Y_1^{(2)}-Y_2^{(2)}\right)s_1^2}{r_1^2} \right)+\frac{Y_2^{(2)}}{4\pi \mu^* (\lambda+2\mu)}\frac{x_1^2 x_2}{ r_1^2 r_2^2}\,.
\notag
\end{align}
It is worth noting that a rigid rotation of the dipoles in Figure \ref{fig:DM} (Right) does not affect the above solution. We emphasize that here as well, the analytical solution \eqref{eq:CDsolution} for a center of dilatation for planar cubic materials is entirely algebraic depending only on the three material parameters $\mu$, $\mu^*$
and $\lambda$.
In the isotropic case, we have that $\mu^*=\mu$ and 
\begin{align}
&\lim_{\mu^* \to \mu} s_m=1 \, , \qquad \lim_{\mu^* \to \mu} {Y_m^{(1)}}=-\lim_{\mu^* \to \mu}{Y_m^{(2)}}=\lambda+\mu \,  , \qquad \lim_{\mu^* \to \mu}q=0 \, , \\
&\lim_{\mu^* \to \mu} r_m=r=\sqrt{x_1^2+x_2^2} \, , \qquad m=1,2\,.
\notag
\end{align}
Accordingly, the solution for the isotropic center of dilation becomes 
\begin{equation}
u_1=\frac{x_1}{4 \pi (\lambda+2\mu)\.r^2} \, , \qquad
u_2=\frac{x_2}{4 \pi (\lambda+2\mu)\.r^2}\, .\label{eq:CDlinear}
\end{equation}
\begin{figure}[h!]
    \centering
    \includegraphics[scale=0.7]{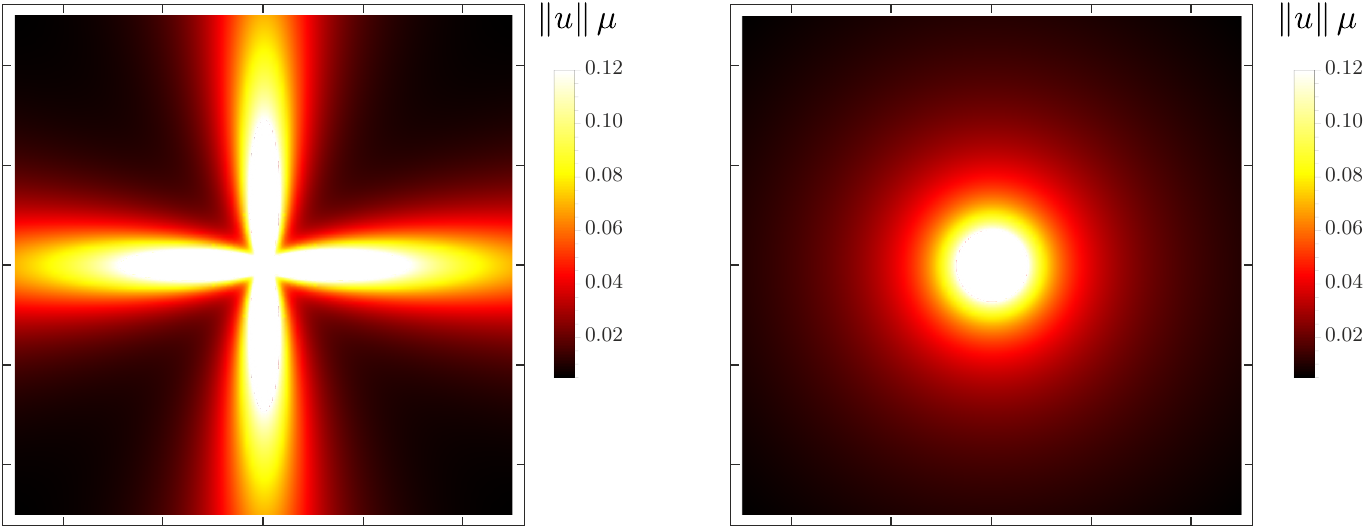}
    \caption{Variation of the normalized modulus of displacement $\norm{u} \mu$ for a center of dilatation (unit strength) (Left) Cubic material ($\mu^*=0.1\. \mu$, $\lambda=0.29\. \mu$), (Right) Isotropic material ($\lambda=0.29\. \mu$).}
    \label{fig:CD}
\end{figure}
%
%
%
%
%
\section{Macroscopic test}
We are now ready to propose a new straightforward approach to finding a best approximating isotropic elasticity tensor. For this, we use the above derived simple isotropic analytical solutions to obtain the two parameters $\mu$ and $\kappa$ of $\Ciso$.
\begin{figure}[h!]
	\pgfmathsetmacro{\r}{1.2}
	\begin{minipage}{0.5\textwidth}
		\centering
		\begin{tikzpicture}
			\draw (-2.5*\r,-2.5*\r) rectangle (2.5*\r,2.5*\r);
			\draw [->] (-1.7*\r,0) -- (1.7*\r, 0) node[below] {$x_1$};
			\draw [->] (0,-1.7*\r) -- (0, 1.7*\r) node[right] {$x_2$};
			\draw [thick](0,0) circle (\r);
			\foreach \a in {0,30,...,360}
				\draw [very thick,->,red] ({\r*cos(\a)},{\r*sin(\a)}) -- ({\r*(cos(\a)-sin(\a))},{\r*(sin(\a)+cos(\a))});
		\end{tikzpicture}	
	\end{minipage}
	\begin{minipage}{0.5\textwidth}
		\centering
		\begin{tikzpicture}
			\draw (-2.5*\r,-2.5*\r) rectangle (2.5*\r,2.5*\r);
			\draw [->] (-1.7*\r,0) -- (1.7*\r, 0) node[below] {$x_1$};
			\draw [->] (0,-1.7*\r) -- (0, 1.7*\r) node[right] {$x_2$};
			\draw[thick] (0,0) circle (\r);
			\foreach \a in {0,30,...,360}
				\draw [very thick,->,red] ({\r*cos(\a)},{\r*sin(\a)}) -- ({1.5*\r*cos(\a)},{1.5*\r*sin(\a)});
		\end{tikzpicture}
	\end{minipage}
	\caption{Instead of applying the load at the origin as done in the analytical solution in Figure \ref{fig:DM}, in the numerical approximation we apply the load on a small circle at the center of the domain. It remains to choose the diameter of the circle in relation to the width of the domain. (Left) Schematic drawing of the statically equivalent reactions in the case of the concentrated couple. (Right) Schematic drawing of the statically equivalent reactions in the case of the center of dilatation.}\label{fig:drawing}
\end{figure}
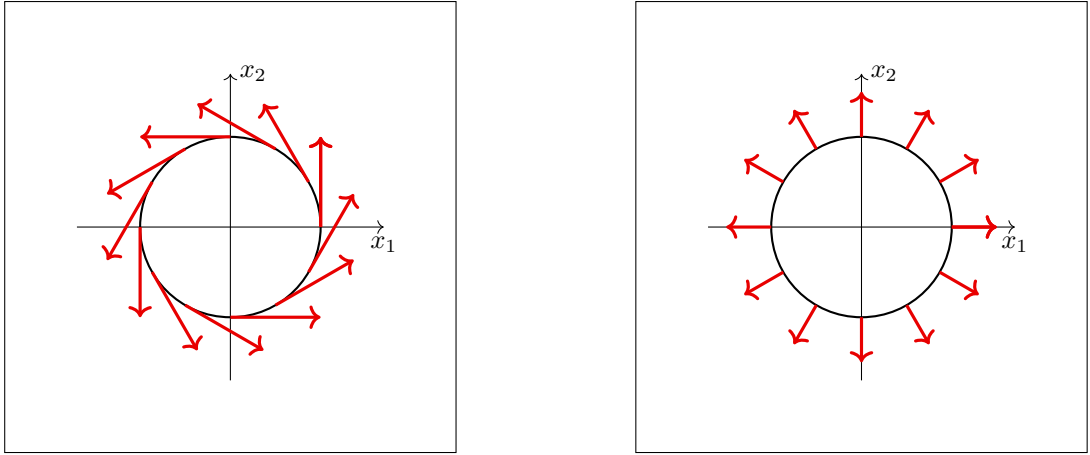
%
%
%
\subsection{Linear elastic solution for the concentrated couple}\label{sec:solutionConcCouple}
We already computed \eqref{eq:CClinear} the linear elastic solution for a concentrated couple. Namely, it is radially symmetric (cf.\ Figure \ref{fig:drawing}) for any given isotropic material and takes the form
\begin{equation}
	u(x_1,x_2)=\frac{1}{4\pi\.\mu\.(x_1^2+x_2^2)}\matr{-x_2\\x_1}.\label{eq:concentratedCouple}
\end{equation}
Thus, the norm of the displacement for any given radius $r=\sqrt{x_1^2+x_2^2}$ is
\begin{equation}
	\norm{u(r)}=\frac{1}{4\pi\.\mu\.(x_1^2+x_2^2)}\sqrt{x_2^2+x_1^2}=\frac{1}{4\pi\.\mu\.r}.\label{eq:concentratedCoupleNorm}
\end{equation}
%
%
%
\subsection{Linear elastic solution for the center of dilatation}\label{sec:centerOfDilatation}
We already computed \eqref{eq:CDlinear} the linear elastic solution for a center of dilatation. It is radially symmetric (cf.\ Figure \ref{fig:drawing}) for any given isotropic material and has the form
\begin{equation}
	u(x_1,x_2)=\frac{1}{4\pi(\mu+\kappa)\.(x_1^2+x_2^2)}\matr{x_1\\x_2}=\frac{1}{4\pi(2\.\mu+\lambda)\.(x_1^2+x_2^2)}\matr{x_1\\x_2},\label{eq:centerOfDilatation}
\end{equation}
since in plane-strain the bulk modulus satisfies $\kappa=\mu+\lambda$, see \cite{agn_gourgiotis2023green}. Thus, the norm of the displacement for any given radius $r=\sqrt{x_1^2+x_2^2}$ is
\begin{equation}
	\norm{u(r)}=\frac{1}{4\pi(\mu+\kappa)(x_1^2+x_2^2)}\sqrt{x_2^2+x_1^2}=\frac{1}{4\pi(\mu+\kappa)\.r}=\frac{1}{4\pi(2\.\mu+\lambda)\.r}.\label{eq:centerOfDilatationNorm}
\end{equation}
%
%
%
\subsection{Geometry of the macroscopic test}
\begin{figure}[h!]
    \centering
	\begin{minipage}{0.45\textwidth}
		\includegraphics[width=\linewidth]{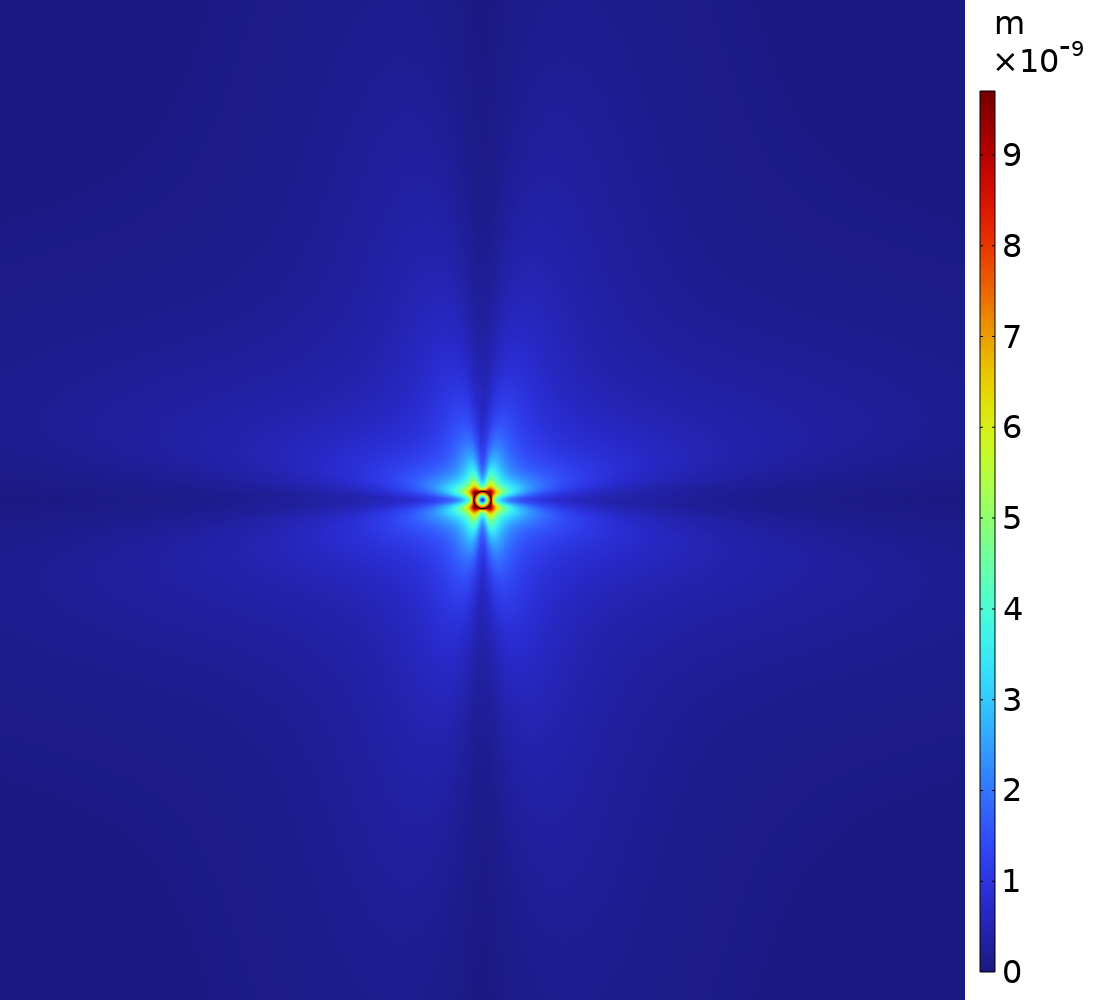}\centering
	\end{minipage}
    \hspace{0.05\textwidth}
	\begin{minipage}{0.45\textwidth}
		\includegraphics[width=\linewidth]{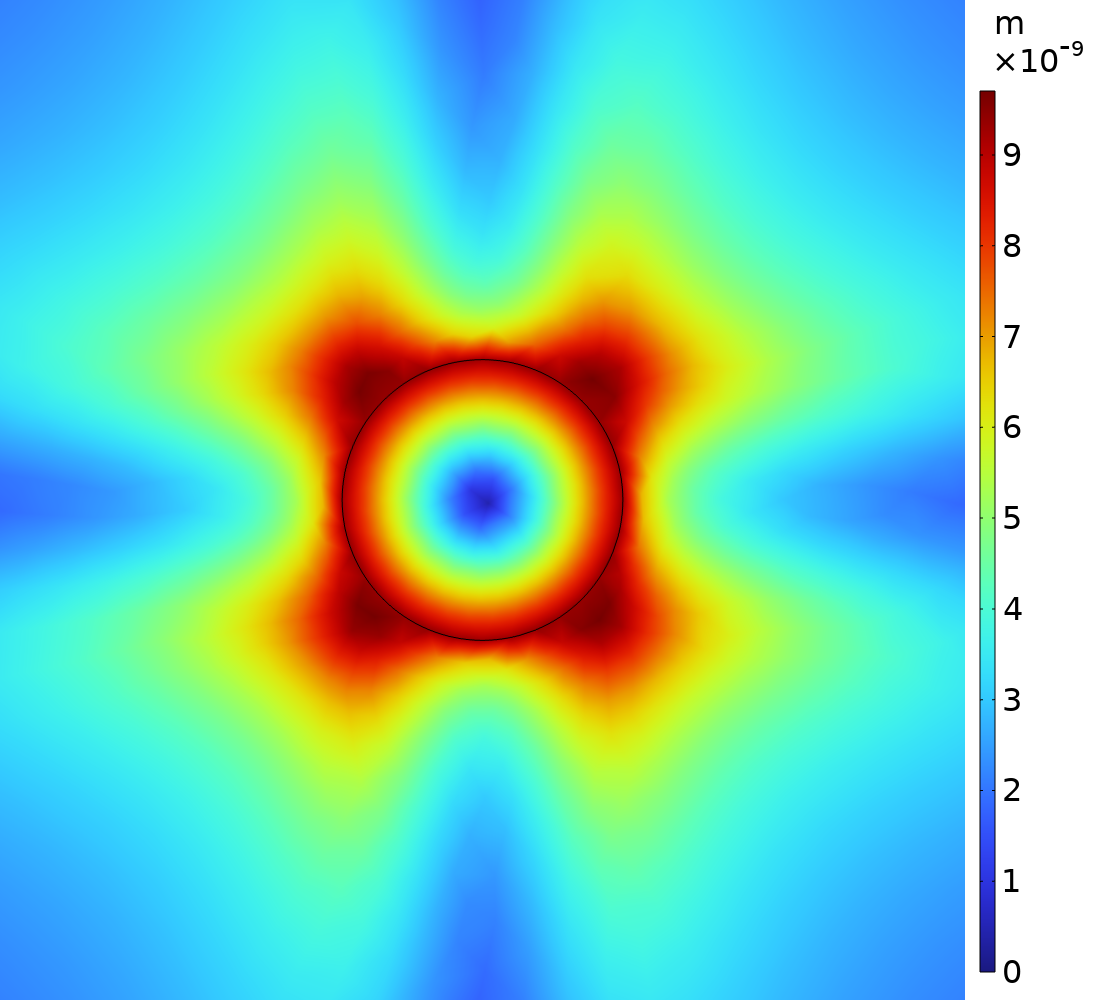}\centering
	\end{minipage}
	\caption{(Left) Displacement for the concentrated couple for the given (anisotropic) Cauchy material $\Caniso$. (Right) Displacement for the concentrated couple for the same specimen zoomed in. The black circle indicates the location where the loads are applied.}\label{fig:concentratedCouple}
\end{figure}

\begin{figure}[h!]
    \centering
	\begin{minipage}{0.45\textwidth}
		\includegraphics[width=\linewidth]{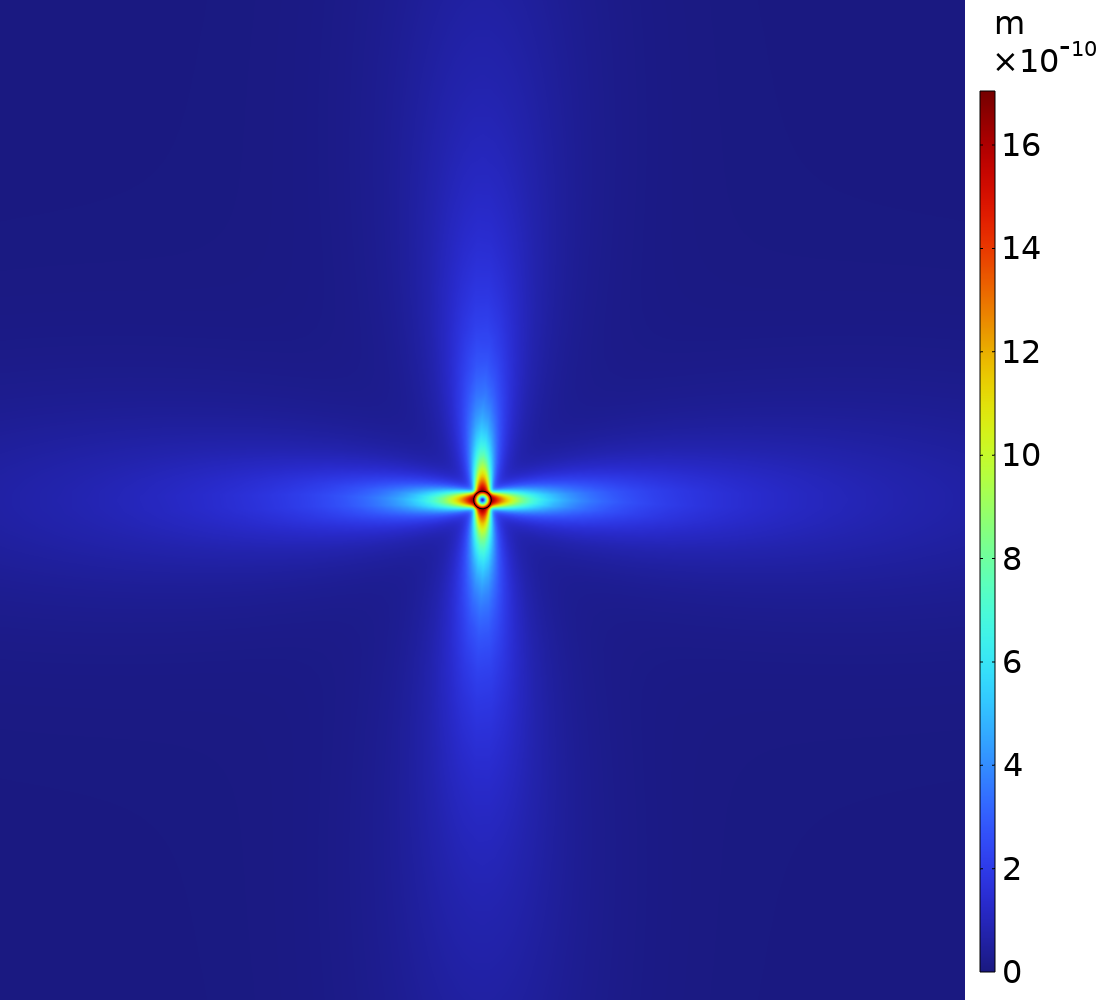}\centering
	\end{minipage}
    \hspace{0.05\textwidth}
	\begin{minipage}{0.45\textwidth}
		\includegraphics[width=\linewidth]{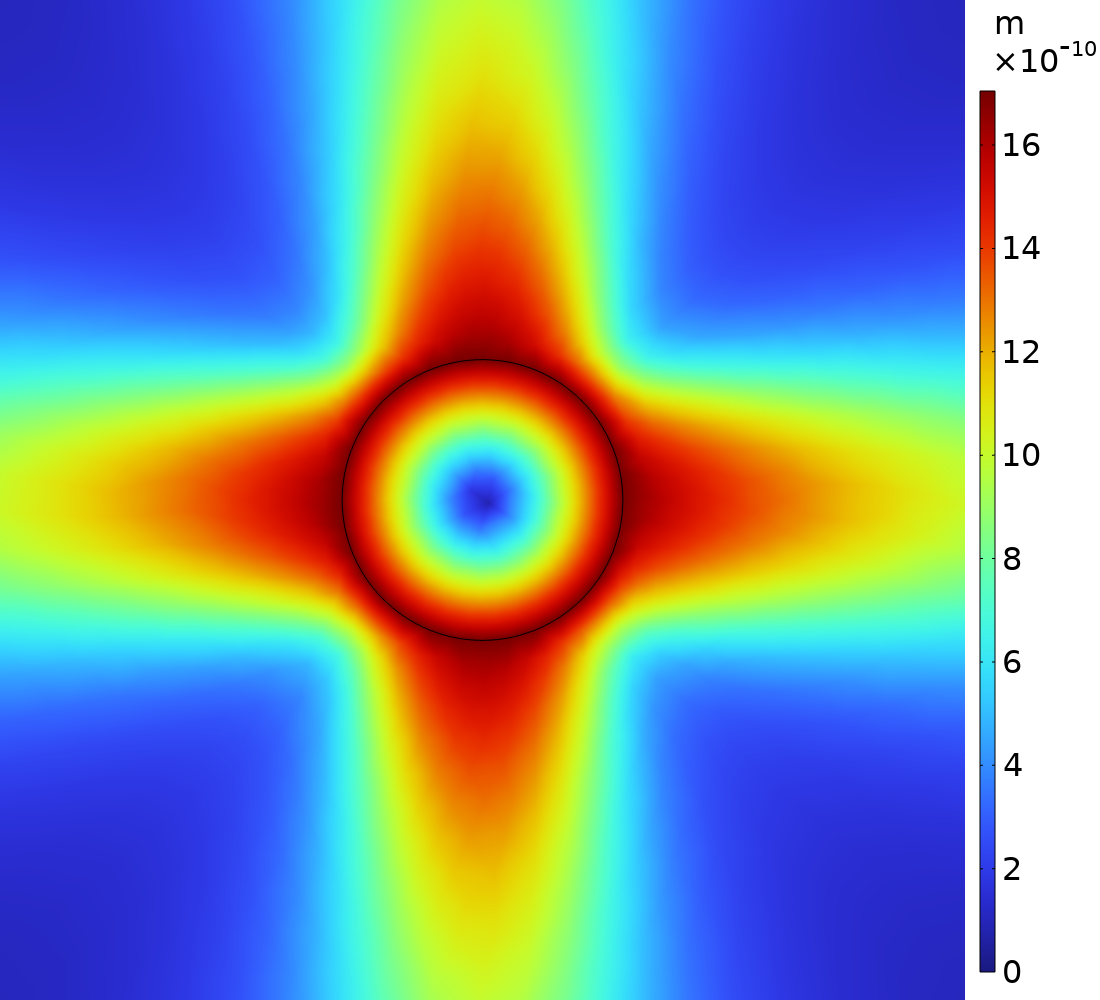}\centering
	\end{minipage}
	\caption{(Left) Displacement for the center of dilatation for the given (anisotropic) Cauchy material $\Caniso$. (Right) Displacement for the center of dilatation of the same specimen zoomed in. The black circle indicates the location where the loads are applied.}\label{fig:centerOfDilatation}
\end{figure}

We consider a homogeneous (but not isotropic) base material described by a given anisotropic elasticity tensor $\Caniso$. While the test works with all symmetry classes, we are currently interested in those with a cubic symmetry class. We perform the two simple static tests (concentrated couple and center of dilation), for which we already gave the isotropic (radially symmetric) solution, cf.\ Section \ref{sec:solutionConcCouple} and \ref{sec:centerOfDilatation}. We also know the analytical solution for the cubic case but we will not use it here to keep the procedure as general as possible.
The corresponding finite element simulation is done in \comsol\;(see Figure \ref{fig:concentratedCouple} and Figure \ref{fig:centerOfDilatation}). Subsequently, a comparison with the analytical displacement solution is used in order to fit the two elastic parameters $\mu$ and $\kappa$ via \mathematica.


For the FEM calculation, without loss of generality, we consider a square finite domain of the length of $1\si\m$. We apply the load with a total force of $1\si\N$ at a sufficiently small diameter of $1\si{\cm}$ in the center of the domain. 
In addition, we specify the corresponding displacement on the outside of the domain to be zero. The validity of the employed ratio of 100:1 between the size of the domain and the diameter of the interior circle stems from and is conformed by extensive numerical testing, cf. Table \ref{table:sizeOfDomain}.
For all FEM calculations, we make use of the automatic physics-controlled mesh from \comsol\;with the extremely fine option ensuring high fidelity of the simulation.

We only know the analytical solution for an isotropic $\Ciso$ (or cubic case) but not for any given (anisotropic) Cauchy material $\Caniso$.
Thus, we are experiencing uncertainties due to the effects of necessary Dirichlet boundary conditions on the outer sides of a finite domain in contrast to an infinite domain with zero Dirichlet boundary conditions at infinity.
Additionally, the best-approximated isotropic displacement using \eqref{eq:concentratedCouple} or \eqref{eq:concentratedCoupleNorm} is unknown as it naturally depends on the very parameters $\muiso$, $\kappaiso$ that we want to determine here with this test. We can also use Norris formula \eqref{eq:norrisLog} in the cubic case but decide against it, in order to keep the method as general as possible, allowing us to find the best approximating isotropic counterpart $\Ciso$ for \textbf{any} anisotropic elasticity tensor $\Caniso$. Thus, we apply zero Dirichlet boundary conditions on all four sides and, precautionarily, only use the displacement data inside a circle with a diameter of $0.5\si{m}$ (half the size of the domain). As per Saint-Venant's principle, this together with a sufficiently large ratio between the size of the domain and the interior circle where the load is applied allows for sufficient accuracy of our method. In section \ref{sec:isotropicCase}, we test different domain sizes in an isotropic setting $\Caniso=\Ciso$ and compare our values $\muiso$, $\kappaiso$ with the original parameters $\mu$ and $\kappa$ from $\Ciso$.
%
%
%
\subsection{Fitting only against the norm $\norm{u(r)}$ of the displacement}\label{sec:fittingNorm}
\begin{figure}[h!]
	\begin{minipage}{0.5\textwidth}
		\includegraphics[width=\linewidth]{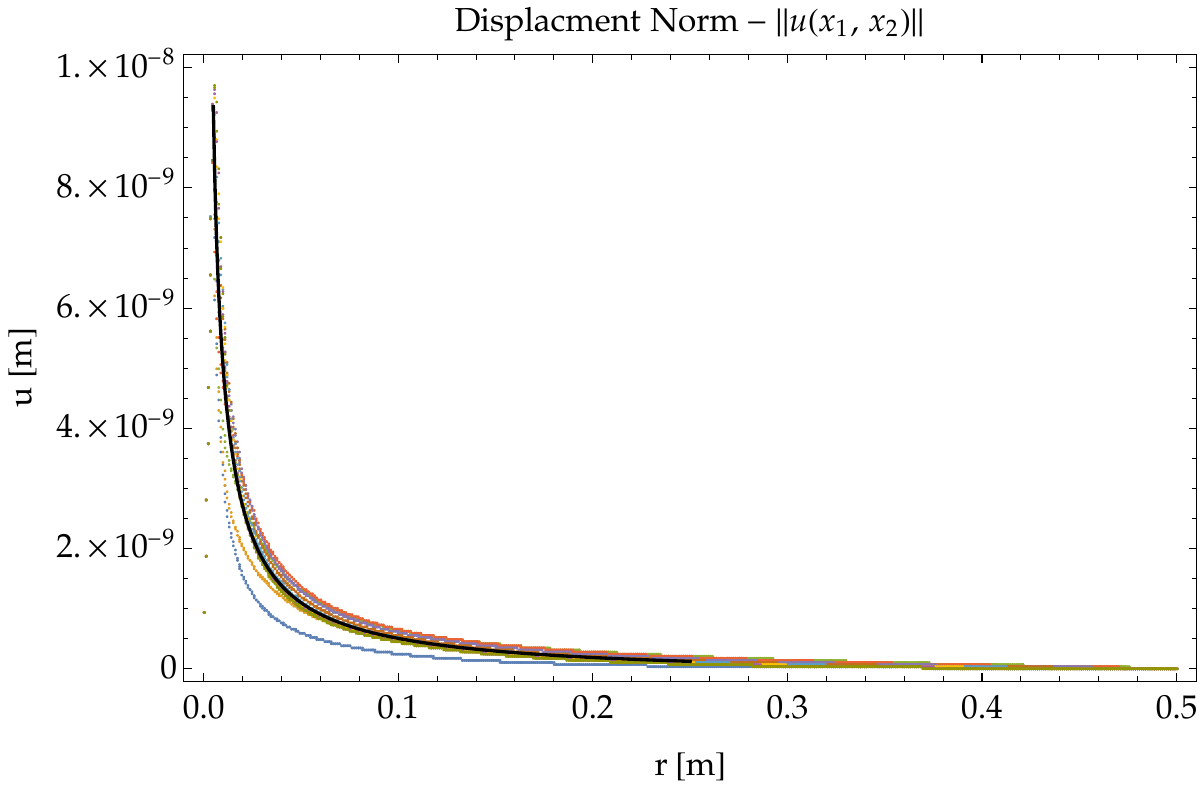}\centering
	\end{minipage}
	\begin{minipage}{0.5\textwidth}
		\includegraphics[width=\linewidth]{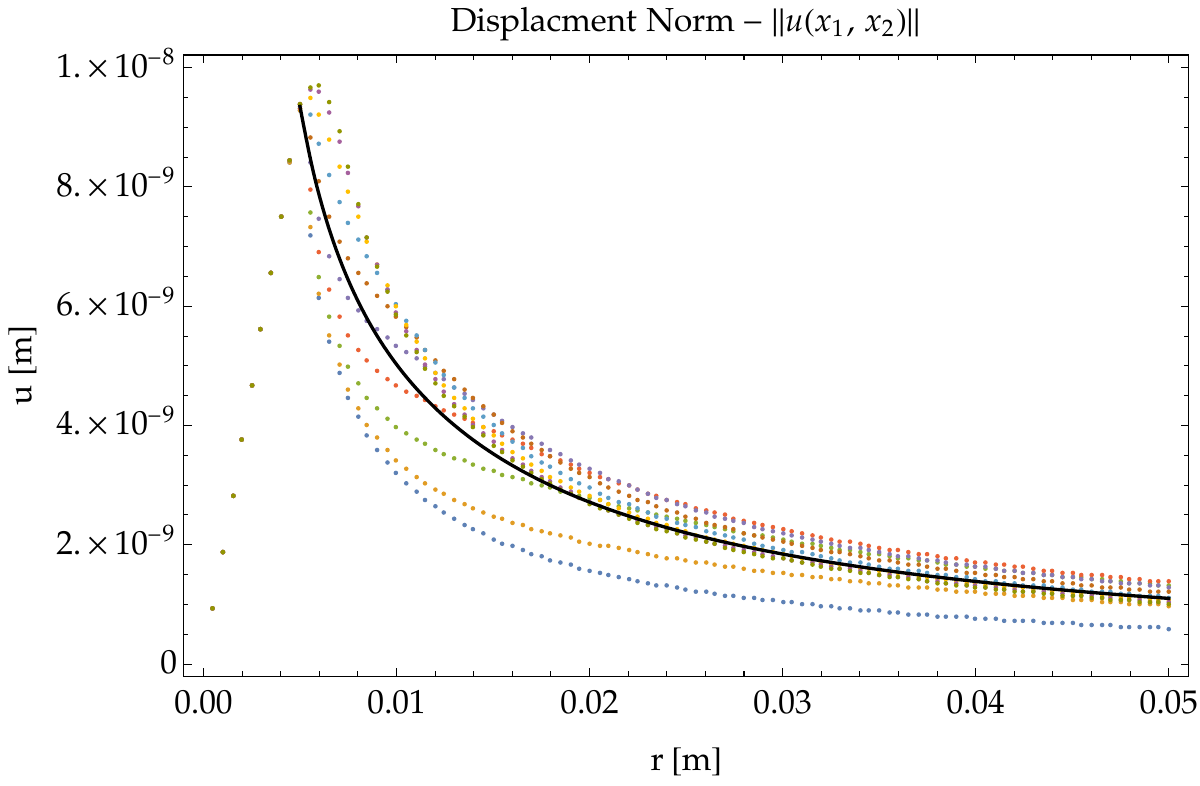}\centering
	\end{minipage}
	\caption{(Left) Norm of the displacement $\norm{u(r)}$ from the FEM calculation for the concentrated couple across different angles $\alpha=0,5,\cdots,45$ (dots in different colors) and its average (black line). (Right) Close-up of the same setting showing only the first 100 data points.}\label{fig:averageAngle}
\end{figure}

We present two different procedures for finding the best approximated isotropic elasticity tensor with quadratic error minimization using \mathematica. In the first one, we only consider the radially averaged norm of the displacement $\norm{u(r)}$. In the second procedure, we take the full displacement solution $u(x_1,x_2)$ into account. For the former, we average the norm of the displacement $\norm{u(r)}$ over all angles\footnote{For the cubic symmetry class considered here, an average between $0^\circ$ and $45^\circ$ would be sufficient.}, cf. Figure \ref{fig:averageAngle}. Note again, that the exact magnitude of the Dirichlet boundary conditions is not known a priori as it depends on the final best approximation isotropic elasticity tensor $\Ciso$. Therefore, we set zero Dirichlet boundary conditions on the outside and only use the data within a circle of half the size of the computed domain for the actual fitting procedure.

First, we find the best approximating isotropic counterpart $\muiso$ by minimizing the mean squared error of the norm of the displacement $\norm{u(r)}$ for the concentrated couple over all points from the anisotropic FEM calculation and the corresponding analytical expression of the isotropic solution \eqref{eq:concentratedCoupleNorm}. With $\muiso$ determined, we obtain the best approximating isotropic counterpart $\kappaiso$ by minimizing the mean squared error of the norm of the displacement $\norm{u(r)}$ for the center of dilatation over all points from the anisotropic FEM calculation and the corresponding analytical expression of the isotropic solution \eqref{eq:centerOfDilatationNorm}. An alternative continuous minimization featuring the logarithmic distance can be found in Appendix \ref{app:geodesics}.
%
%
%
%
\subsection{Fitting against the full  displacement solution $u(x_1,x_2)$}\label{sec:fittingAll}
In our second approach, we minimize the difference between the full displacement field $u(x_1,x_2)$ of its anisotropic FEM calculations and the corresponding analytical expressions. In this scenario, the FEM calculation itself remains unchanged but a different set of data points is used. While polar coordinates are used in the first approach, here, we consider a regular (Cartesian) grid instead. This creates evenly distributed data points over the whole domain, in contrast to a clear accumulation of points in the proximity of the origin when using polar coordinates. The second main difference is the absence of averaging over the angle\footnote{Depending on the symmetry class of $\Caniso$ it is still possible to reduce the total amount of points considered, e.g. for the cubic symmetry class considered here, a section between $0^\circ$ and $45^\circ$ is still sufficient.}. Again, we only use data points inside a circle of half the size of the computed domain for the actual fitting procedure to avoid boundary effects due to imperfect values chosen (zero Dirichlet boundary conditions on the outside).

We start by finding the best approximating isotropic counterpart $\muiso$ via minimization of the mean squared error of the displacement for the concentrated couple, using both displacement components $u_1(x_1,x_2)$ and $u_2(x_1,x_2)$ separately. The minimization is done across all points from the anisotropic FEM calculation and the corresponding analytical expression for the isotropic solution \eqref{eq:concentratedCouple}. Next, we find the best approximating isotropic counterpart $\kappaiso$ by minimizing the mean squared error of the norm of the displacement for the center of dilatation, where we again employ both components $u_1(x_1,x_2)$ and $u_2(x_1,x_2)$ separately. This is done using $\muiso$ from the former step and a minimization over all points from the anisotropic FEM calculation and the corresponding analytical expression for the isotropic solution \eqref{eq:centerOfDilatation}.
%
%
%
\section{Results}
For most of the calculations and all the plots (except the following isotropic consistency check), we consider the anisotropic elasticity tensor $\Caniso$ with cubic symmetry \eqref{eq:voigtNotation} and the following material parameters \cite{agn_sarhil2023identification,agm_sarhil2023size,agn_barbagallo2017transparent}, approximating the macroscopic material of metamaterial with cubic symmetry
\begin{itemize}
	\item the bulk modulus \hspace{1.4cm}$\kappa=7.645\,[\si{\GPa}]$,
	\item the first shear modulus \hspace{0.55cm}$\mu=5.901\,[\si{\GPa}]$,
	\item the second shear modulus $\mus=0.626\,[\si{\GPa}]$.
\end{itemize}
While this specific set of parameters shows a close similarity in values for $\muiso$ and $\kappaiso$ in our two fitting procedures and the analytical values from Norris, we also test other anisotropic elasticity tensors $\Caniso$ where differences between all three methods are clearly visible.
We emphasize again that for the cubic case only, we could omit the numerical FEM approach and use the derived analytical expressions instead, cf.\ Section \ref{sec:GreenFunction}. While we fitted the best approximation isotropic counterpart $\Ciso$ only using the FEM calculation, we validated the displacements with the analytical expression of the concentrated couple \eqref{eq:CCsolution} and the center of dilatation \eqref{eq:CDsolution}, respectively.

%
%
%
\subsection{The isotropic case - consistency check}\label{sec:isotropicCase}
We start with the isotropic case, i.e. $\mu^*=\mu$ as the simple example where the isotropic exact solution is known a priori. Given the analytical solutions \eqref{eq:concentratedCouple} and \eqref{eq:centerOfDilatation}, we can compute the correct displacement $u(x_1,x_2)$ at the boundary of any domain.
Thus, we can prescribe the exact value at the outside boundary depending on $\muiso$ instead of zero Dirichlet boundary conditions, matching the force applied inside.\footnote{The same consistency check can be done in the cubic case using the analytical solutions \eqref{eq:CCsolution} and \eqref{eq:CDsolution}.} In general, for an anisotropic tensor $\C_{\rm aniso}$ this value is not known a priori because it depends on the isotropic best approximated $\muiso$ which is to be determined. 

While for smaller domains it seems to be necessary to adjust the value of the boundary displacement correctly (which can always be done recursively), we circumvent the problem by increasing the considered domain. Therein, we change the overall size but not the radii of the applied force nor the maximum distance of points used for the fitting procedure. This allows us to reduce the effect of inaccuracies at the boundary as much as required which is attributed to Saint-Venant's principle. Consequently, for a sufficiently large domain, we can set zero Dirichlet displacements at the outside boundary.

\begin{figure}[h!]
    \centering
	\begin{minipage}{0.45\textwidth}
		\includegraphics[width=\linewidth]{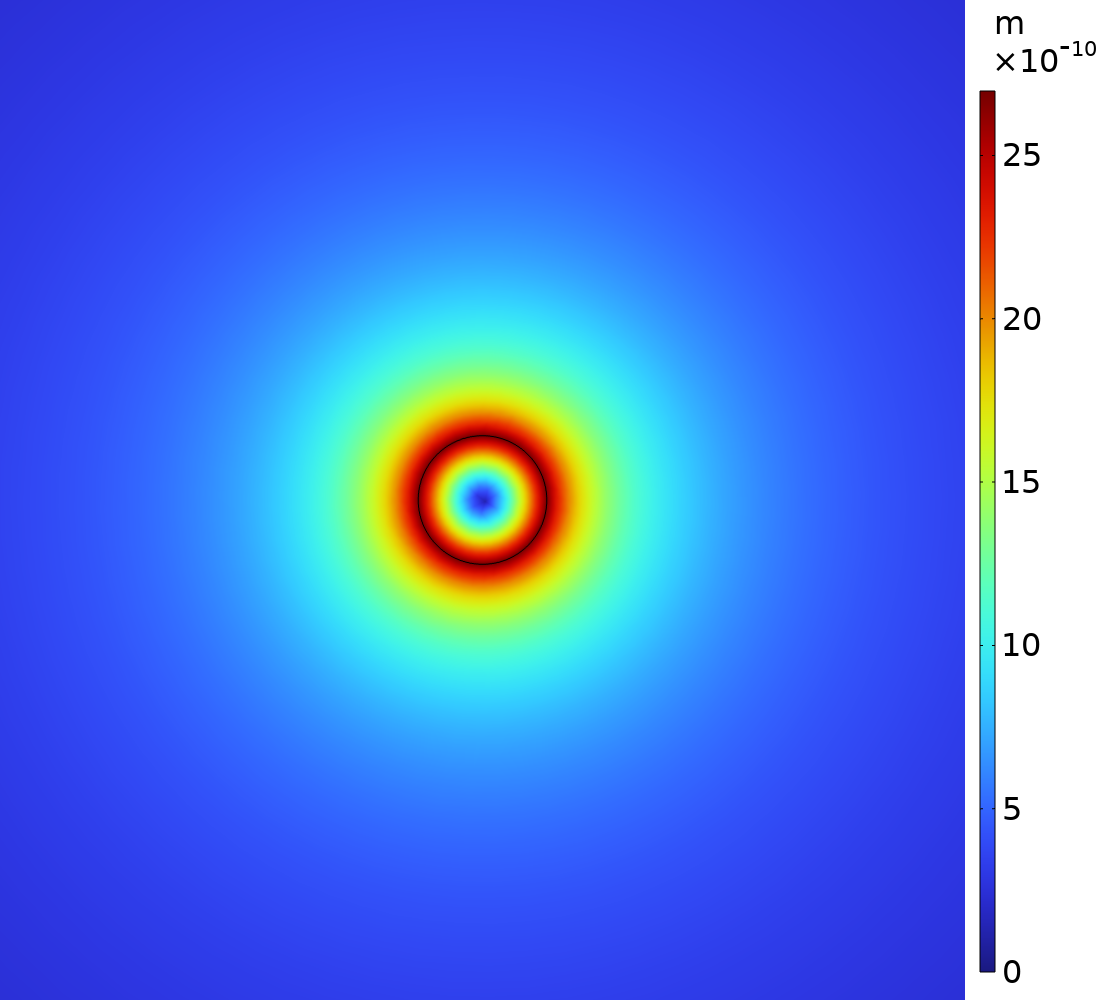}
		\caption{Displacement for the concentrated couple in the isotropic case $\mu^*=\mu$, zoomed in. The black circle indicates the location where the loads are applied.}
	\end{minipage}
	\hfill
	\begin{minipage}{0.5\textwidth}
		\centering
		\begin{tabularx}{\textwidth}{c|llll}
			size of domain & $0.5\si{m}$ & \color{green}$1\si{m}$ & $2.5\si{m}$ & $10\si{m}$\\[0.2em]\hline\\[-0.7em]
			$\displaystyle\frac{\muiso}{\mu}$ (norm) & 1.0174 & 1.0042 & 1.0006 & 1.0001 \\[0.2em]\hline\\[-0.7em]
			$\displaystyle\frac{\kappaiso}{\kappa}$ (norm) & 1.0156 & 1.0039 & 1.0006 & 1.00001 \\[0.2em]\hline\\[-0.7em]
			$\displaystyle\frac{\muiso}{\mu}$ (disp) & 1.1263 & 1.0289 & 1.0049 & 1.0003 \\[0.2em]\hline\\[-0.7em]
			$\displaystyle\frac{\kappaiso}{\kappa}$ (disp) & 1.1149 & 1.0265 & 1.0038 & 1.0002
		\end{tabularx}
		\captionof{table}{Different sizes of domains are considered while the diameter of the interior circle is fixed at $1\si{cm}$ and only points up to the maximum distance of $0.5\si{m}$ are used for the fitting. 
		These values are computed by minimizing the norm of the displacement $\norm{u(r)}$ (norm) and the full displacement field $u(x_1,x_2)$ (disp), respectively, and show the discrepancy to the analytical solution of an infinite domain.}\label{table:sizeOfDomain}
	\end{minipage}
\end{figure}
In Table \ref{table:sizeOfDomain}, different sizes of domains are considered and the numerical solution of both of our methods, i.e. considering the norm of the displacement $\norm{u(r)}$, or the displacement field $u(x_1,x_2)$ itself, are compared with the analytical solution for $\muiso=\mu$ and $\kappaiso=\kappa$. We observe that the error decreases significantly for larger domains such that for our calculations, we choose a domain size of $1\si{m}$, i.e. a ratio of 100:1 to the diameter of the interior circle where the loads are applied.
%
%
%
\subsection{The best approximating $\Ciso$ for the norm of the displacement $\norm{u(r)}$}
\begin{figure}[h!]
	\begin{minipage}{0.5\textwidth}
		\includegraphics[width=\linewidth]{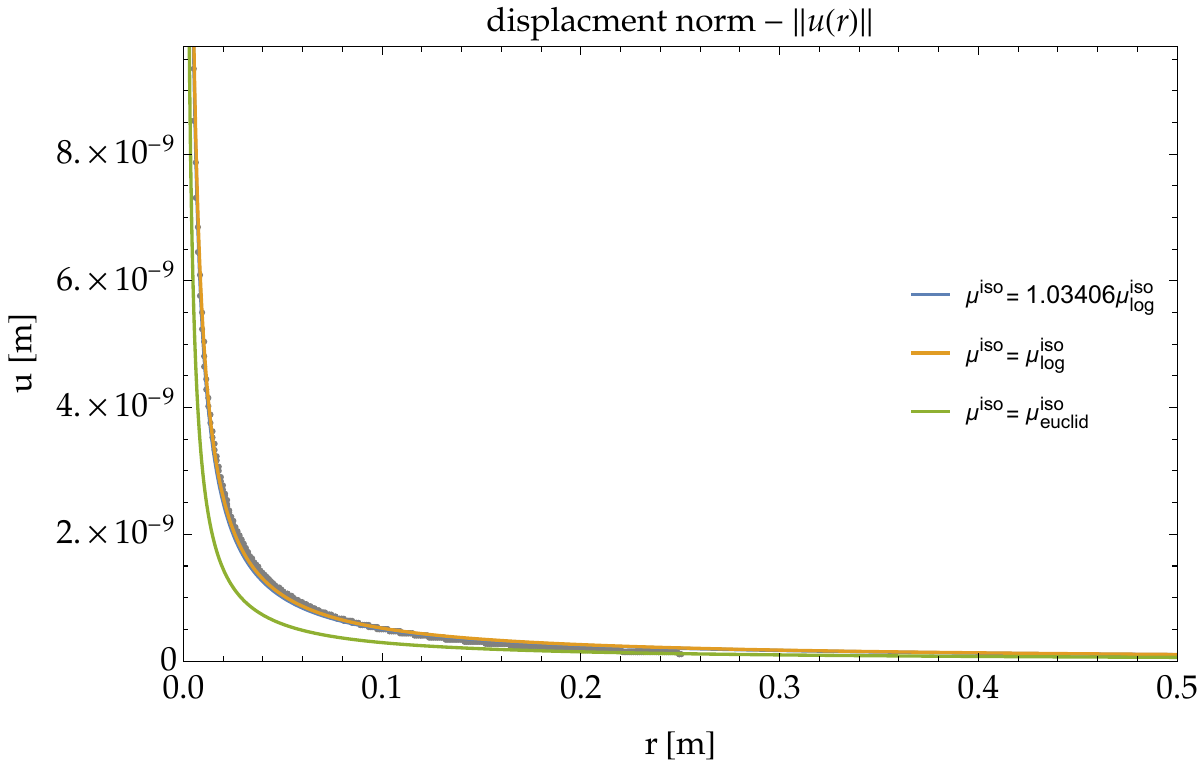}\centering
	\end{minipage}
	\begin{minipage}{0.5\textwidth}
		\includegraphics[width=\linewidth]{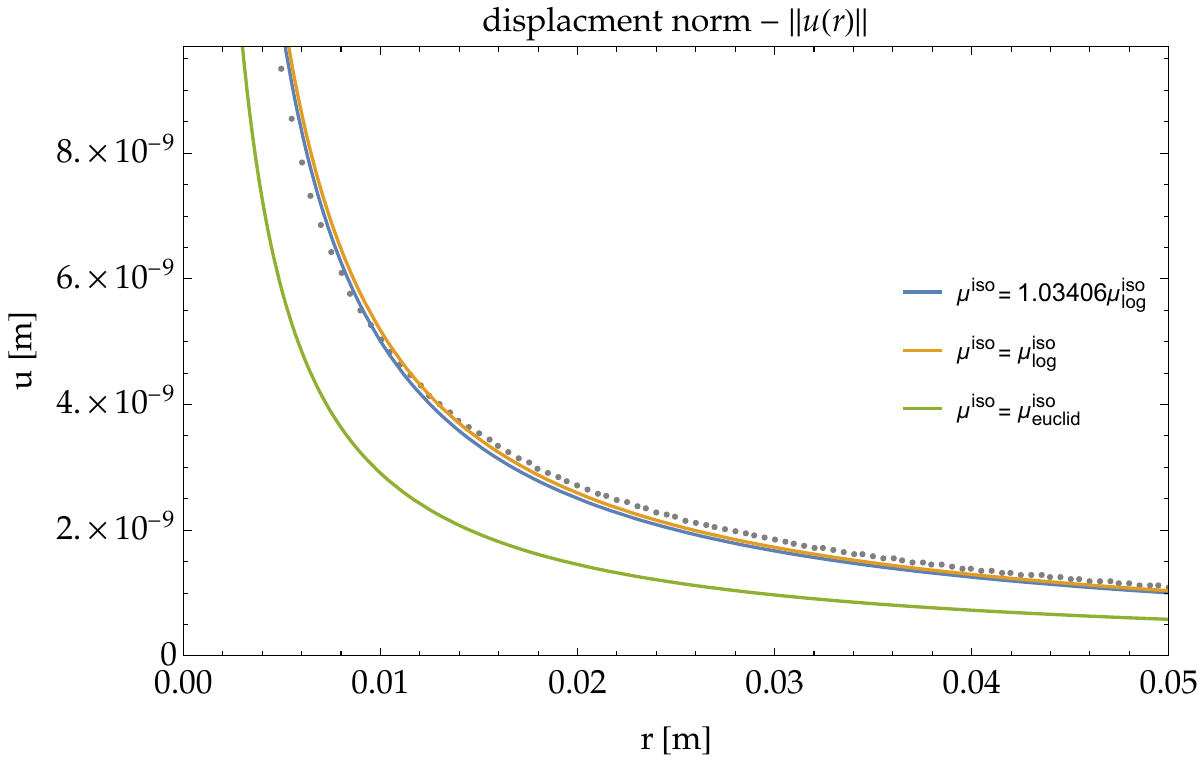}\centering
	\end{minipage}
	\caption{Comparison for the concentrated couple between our numerical data averaged over the angle (black dots), our best isotropic approximation (blue line), and the two different analytical values $\mulogiso$ (orange line) and $\mueuclidiso$ (green line) from Norris. (Left) The plot shows values for the norm of the displacement $\norm{u(r)}$ for the whole domain (only the first half was used for the fitting). (Right) A close-up of the same plot.}\label{fig:fitConcentratedCouple}
\end{figure}

\begin{figure}[h!]
	\begin{minipage}{0.5\textwidth}
		\includegraphics[width=\linewidth]{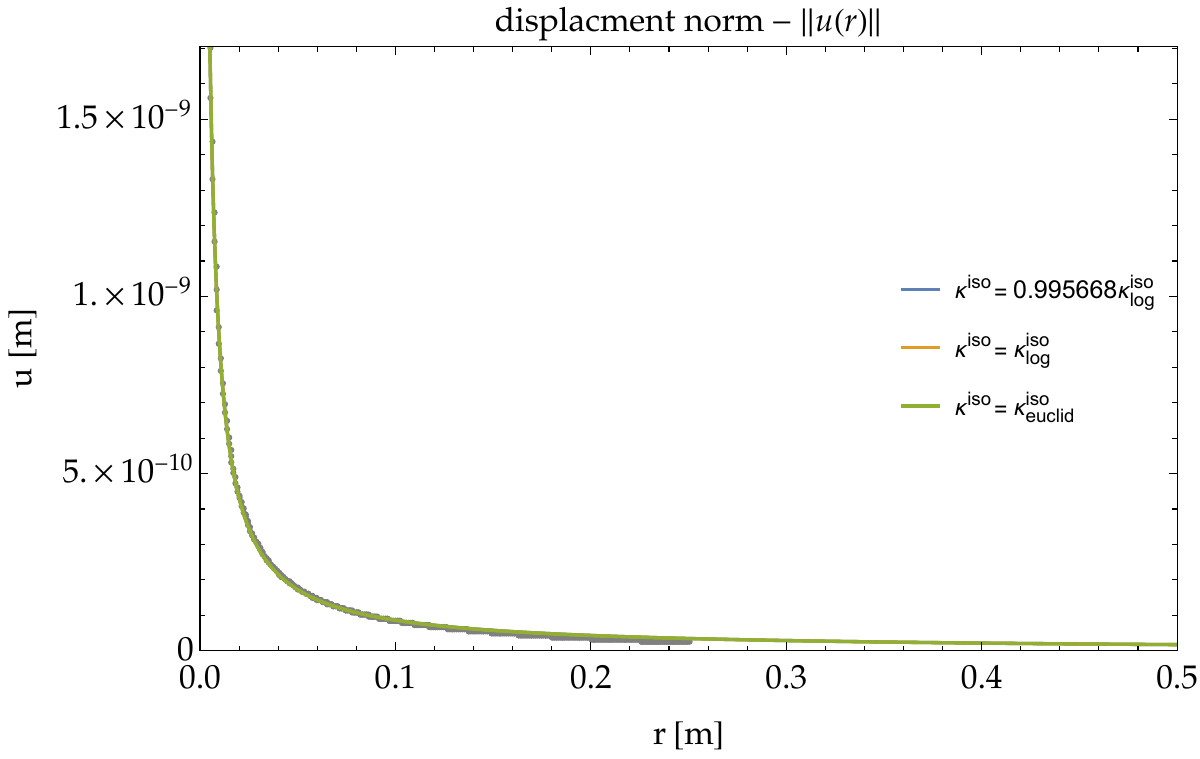}\centering
	\end{minipage}
	\begin{minipage}{0.5\textwidth}
		\includegraphics[width=\linewidth]{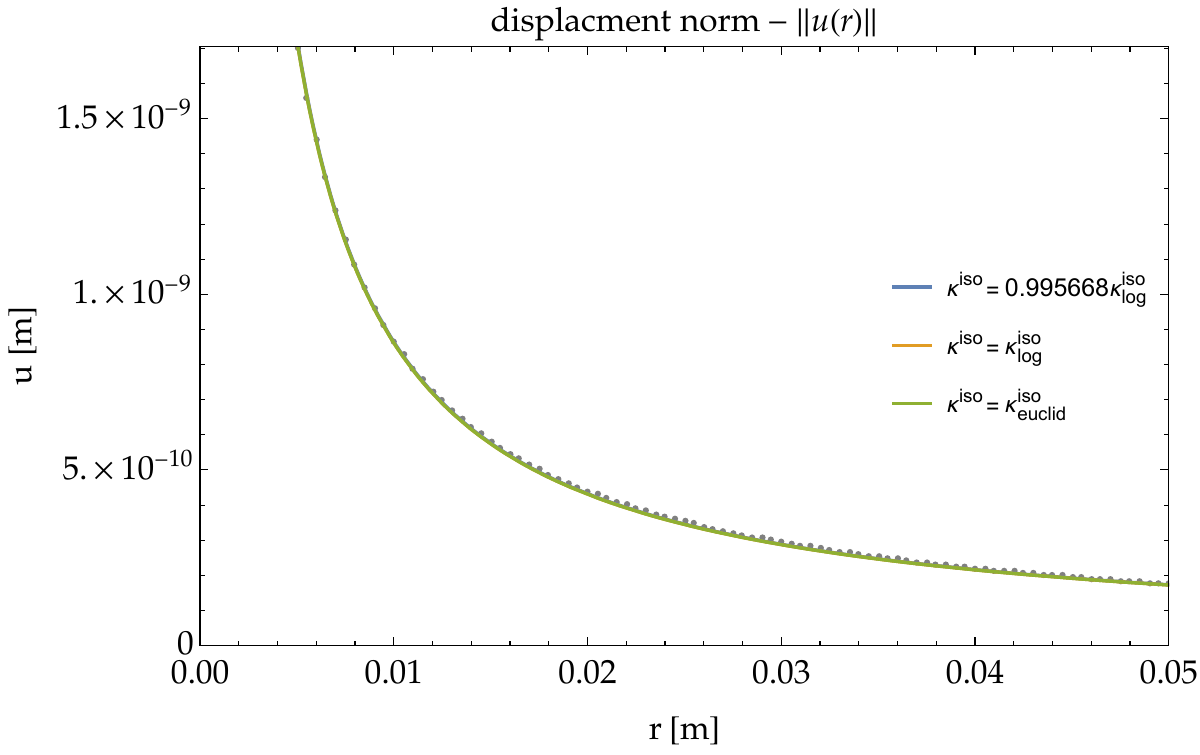}\centering
	\end{minipage}
	\caption{Comparison for the center of dilatation between our numerical data averaged over the angle (black dots), our best isotropic approximation (blue line), and the two different analytical values $\kappalogiso$ (orange line) and $\kappaeuclidiso$ (green line) from Norris. Note that $\kappalogiso=\kappaeuclidiso=\kappa$. (Left) The plot shows values for the norm of the displacement $\norm{u(r)}$ for the whole domain (only the first half was used for the fitting). (Right) A close-up of the same plot.}\label{fig:fitCenterOfDilatation}
\end{figure}

Following our method described in Section \ref{sec:fittingNorm}, we start with the anisotropic elasticity tensor $\Caniso$ introduced above to find the best approximating isotropic shear modules $\muiso$ from the concentrated couple (see Figure \ref{fig:fitConcentratedCouple}), as well as the best approximating isotropic bulk modules $\kappaiso$ from the center of dilatation (see Figure \ref{fig:fitCenterOfDilatation}).

Although the considered anisotropic tensor $\Caniso(\kappa,\mu,\mus)$ is highly non-isotropic with a ratio of almost 10 between $\mu$ and $\mus$, we note that our best approximating isotropic tensor $\Ciso(\kappaiso,\muiso)$ is very close to the logarithmic solution of Norris $\Ciso(\kappa,\mulogiso)$. In particular, the bulk-modulus $\kappa$ remains almost unchanged. In Table \ref{tab:resultNorm} we repeat the fitting for different anisotropic tensors $\Caniso$ by changing $\mus$ while keeping the values for $\kappa$ and $\mu$. Again, the best approximating bulk-modulus $\kappaiso$ is almost identical to the original bulk-modulus $\kappa$ for all ratios considered. However, the computed best approximating shear modulus $\muiso$ differs greatly from the logarithmic solution $\mulogiso$ if $\mus>\mu$.

\begin{table}[h!]
	\centering
	\begin{tabularx}{12cm}{c|lllll}
		$\displaystyle\mus$ & \color{green}$0.626\,[\si{\GPa}]$ & $1.967\,[\si{\GPa}]$ & $5.901\,[\si{\GPa}]$ & $17.70\,[\si{\GPa}]$ & $59.01\,[\si{\GPa}]$\\[0.2em]\hline\\[-0.7em]
		$\displaystyle\frac{\mus}{\mu}$ & 0.106 & $\displaystyle\frac13$ & 1 & 3 & 10 \\[0.2em]\hline\\[-0.7em]
		$\displaystyle\frac{\muiso}{\mulogiso}$ & 1.03406 & 1.04654 & 1.0042 & 0.808581 & 0.494873 \\[0.2em]\hline\\[-0.7em]
		$\displaystyle\frac{\kappaiso}{\kappa}$  & 0.99568 & 1.00434 & 1.00386 & 1.00302 & 1.00628
	\end{tabularx}
	\caption{Values of the best approximating $\Ciso$ minimizing $\norm{u(r)}$ using our method for different $\Caniso$ where we fix the bulk modulus $\kappa=7.645\,[\si{\GPa}]$ and the first shear modulus $\mu=5.901\,[\si{\GPa}]$ but consider different values for $\mus$. The first row marked in green corresponds to the values of a homogenized macroscopic metamaterial described in \cite{agn_sarhil2023identification,agn_barbagallo2017transparent}}\label{tab:resultNorm}
\end{table}
%
%
%
\subsection{The best approximating $\Ciso$ for the displacement $u(x_1,x_2)$}
\begin{figure}[h!]
	\begin{minipage}{0.45\textwidth}
		\includegraphics[width=\linewidth]{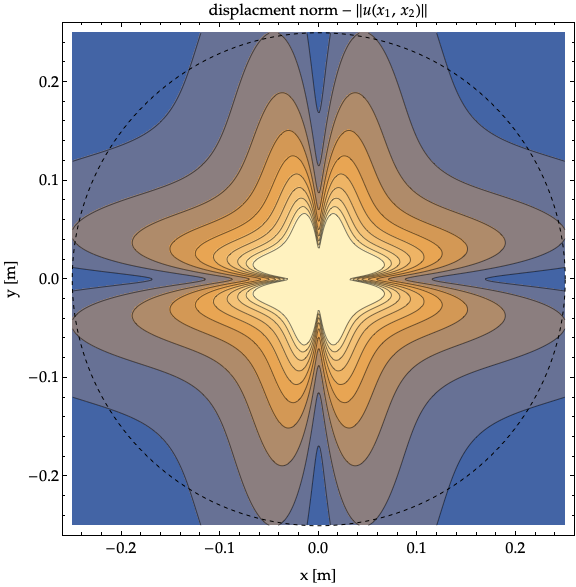}\centering
	\end{minipage}
	\begin{minipage}{0.45\textwidth}
		\includegraphics[width=\linewidth]{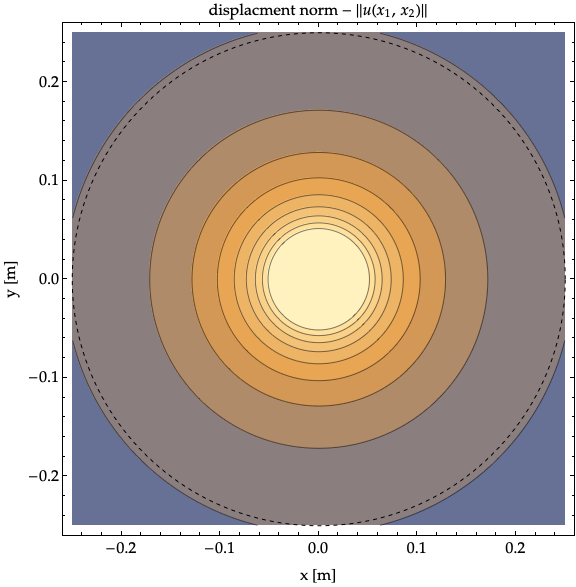}\centering
	\end{minipage}
	\begin{minipage}{0.09\textwidth}
		\includegraphics[width=\linewidth]{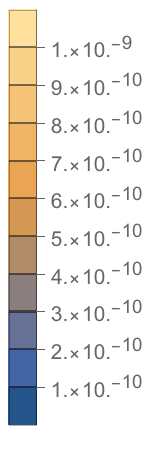}\centering
	\end{minipage}
	\caption{Visualization of the concentrated couple with the color indicating the amplitude of the norm of the displacement $\norm{u(x_1,x_2)}$, only showing the circular area inside the domain which is used for the fitting procedure. (Left) The anisotropic behavior is computed with FEM. (Right) The best approximating isotropic behavior with quadratic error minimization using \mathematica.}\label{fig:fitConcentratedCoupleAll}
\end{figure}

\begin{figure}[h!]
	\begin{minipage}{0.45\textwidth}
		\includegraphics[width=\linewidth]{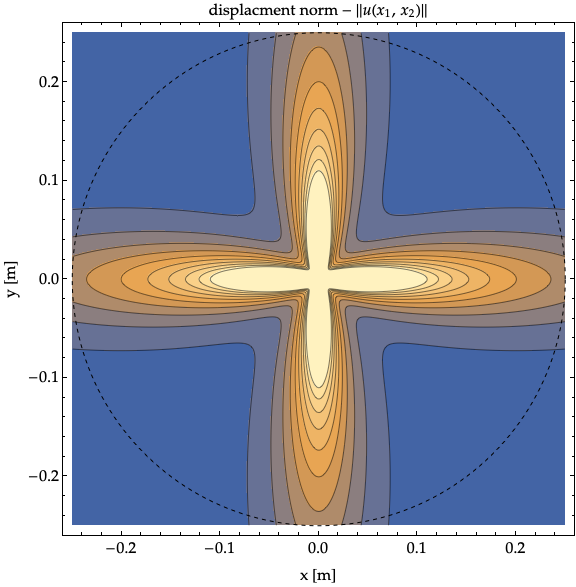}\centering
	\end{minipage}
	\begin{minipage}{0.45\textwidth}
		\includegraphics[width=\linewidth]{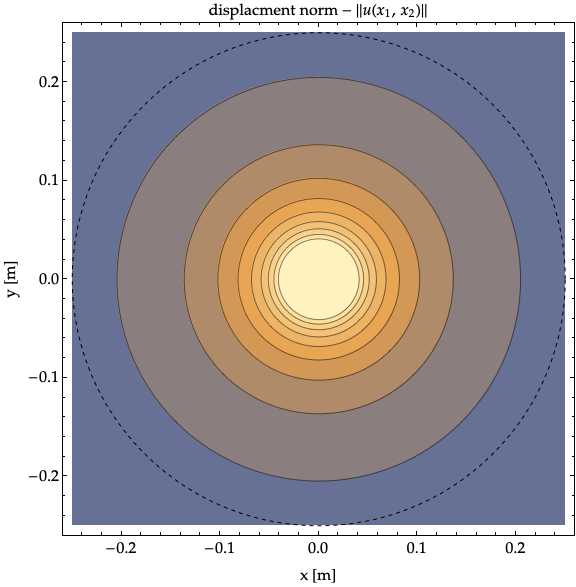}\centering
	\end{minipage}
	\begin{minipage}{0.09\textwidth}
		\includegraphics[width=\linewidth]{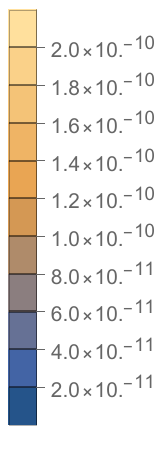}\centering
	\end{minipage}
	\caption{Visualization of the center of dilatation with the color indicating the amplitude of the norm of the displacement $\norm{u(x_1,x_2)}$, only showing the circular area inside the domain which is used for the fitting procedure. (Left) The anisotropic behavior is computed with FEM. (Right) The best approximating isotropic behavior with quadratic error minimization using \mathematica.}\label{fig:fitCenterOfDilatationAll}
\end{figure}

Following our method described in Section \ref{sec:fittingAll}, we start with the anisotropic elasticity tensor $\Caniso$ introduced above to find the best approximating isotropic shear modules $\muiso$ from the concentrated couple (see Figure \ref{fig:fitConcentratedCoupleAll}), as well as the best approximating isotropic bulk modules $\kappaiso$ from the center of dilatation (see Figure \ref{fig:fitCenterOfDilatationAll}).

Here, we derive a best approximating isotropic tensor $\Ciso(\kappaiso,\muiso)$ with the values of $\muiso=1.24043\.\mulogiso$  and $\kappaiso=1.03999\.\kappa$. Thus, our solution clearly differs from the logarithmic solution of Norris and also our previous method of minimizing the norm of the displacement $\norm{u(r)}$. Note that the best approximating bulk-modulus $\kappaiso$ only changes slightly, partly due to boundary effects, cf.\ Table \ref{table:sizeOfDomain}. In Table \ref{tab:resultDisplacement} we repeat the fitting for different anisotropic tenors $\Caniso$ by changing $\mus$ while keeping the values for $\kappa$ and $\mu$ fixed. Again, the best approximating bulk-modulus $\kappaiso$ remains almost unchanged for all ratios considered. However, the computed best approximating shear modulus $\muiso$ differs greatly from Norris' logarithmic solution $\mulogiso$ if $\mus>\mu$.

\begin{table}[h!]
	\centering
	\begin{tabularx}{12cm}{c|lllll}
		$\displaystyle\mus$ & \color{green}$0.626\,[\si{\GPa}]$ & $1.967\,[\si{\GPa}]$ & $5.901\,[\si{\GPa}]$ & $17.70\,[\si{\GPa}]$ & $59.01\,[\si{\GPa}]$\\[0.2em]\hline\\[-0.7em]
		$\displaystyle\frac{\mus}{\mu}$ & 0.106 & $\displaystyle\frac13$ & 1 & 3 & 10 \\[0.2em]\hline\\[-0.7em]
		$\displaystyle\frac{\muiso}{\mulogiso}$ & 1.24043 & 1.11925 & 1.0289 & 0.863217 & 0.575421 \\[0.2em]\hline\\[-0.7em]
		$\displaystyle\frac{\kappaiso}{\kappa}$  & 1.03999 & 1.03357 & 1.02649 & 1.01104 & 1.03283
	\end{tabularx}
	\caption{Values of the best approximating $\Ciso$, minimizing $u(x_1,x_2)$ using our method for different $\Caniso$, where we fix the bulk modulus $\kappa=7.645\,[\si{\GPa}]$ and the first shear modulus $\mus=0.626\,[\si{\GPa}]$ but consider different values for $\mus$. The first row marked in green corresponds to the values of a homogenized macroscopic metamaterial described in \cite{agn_sarhil2023identification,agn_barbagallo2017transparent}}\label{tab:resultDisplacement}
\end{table}
%
%
%
\section{Discussion and application}\label{sec:discussion}
\begin{figure}[h!]
	\begin{minipage}{0.3\textwidth}
		\includegraphics[width=\linewidth]{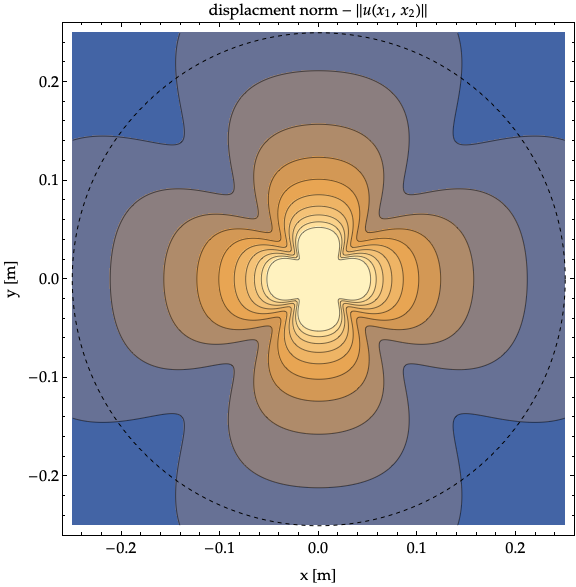}\centering
	\end{minipage}
	\begin{minipage}{0.3\textwidth}
		\includegraphics[width=\linewidth]{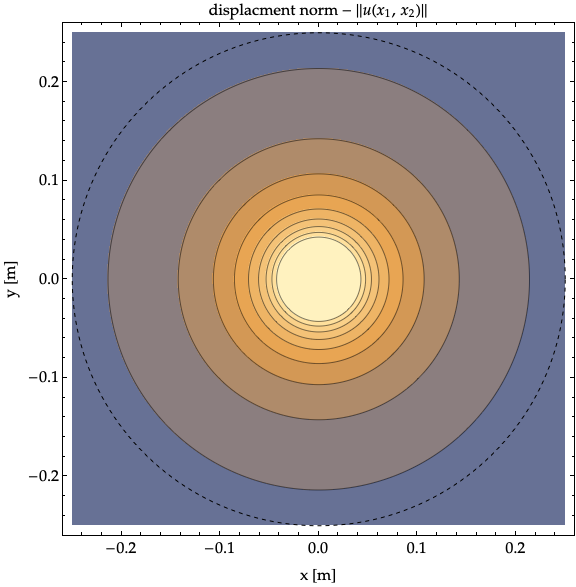}\centering
	\end{minipage}
	\begin{minipage}{0.3\textwidth}
		\includegraphics[width=\linewidth]{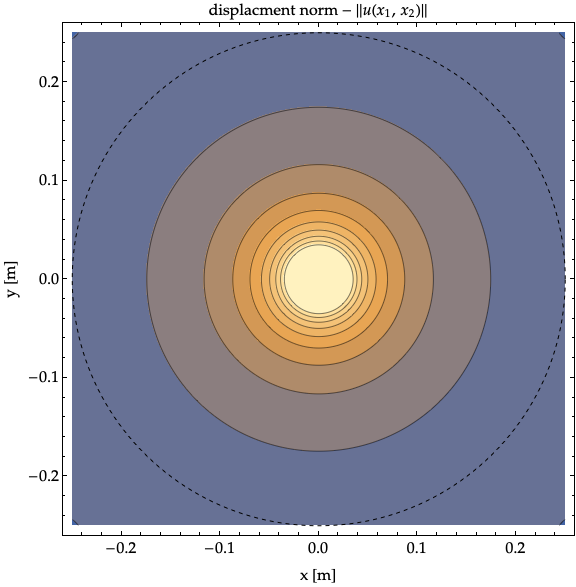}\centering
	\end{minipage}
	\begin{minipage}{0.08\textwidth}
		\includegraphics[width=\linewidth]{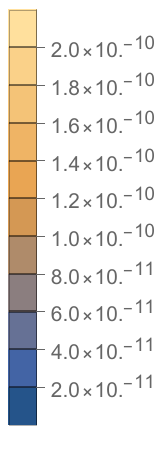}\centering
	\end{minipage}
	\caption{Comparison between the original (anisotropic) material behavior with $\mus=3\mu$ under a concentrated couple (left), our best approximating isotropic behavior (middle) and the isotropic behavior by Norris minimizing the logarithmic distance.}\label{fig:Comparsion}
\end{figure}

We conclude that one should be extremely careful when replacing an anisotropic elasticity tensor $\Caniso$ with an approximating isotropic counterpart $\Ciso$.
There are many different methods available \cite{cavallini1999best,norris2006isotropic,moakher2006closest,kochetov2009estimating,antonelli2022distance,azzi2023distance} and it isn't easy to compare different best-approximating isotropic tensors, cf.\ Figure \ref{fig:Comparsion}.

In the main setting used in this paper, i.e. cubic planar linear elasticity, we already have an easy-to-compute best-approximating isotropic tensor with the Norris formula a \eqref{eq:norrisLog}.

However, we want to emphasize two main advantages of our method for finding $\Ciso$:
\begin{itemize}
	\item[1.)] We have a coherent method of using two simple material tests, the concentrated couple and the center of dilatation, to identify the two physical parameters $\mu$ and $\kappa$ which can be easily visualized.
	\item[2.)] Our method is applicable without modification to any class of symmetry, i.e.\ we can allow arbitrary anisotropic elasticity tensors $\Caniso\in\Symp(6)$.
\end{itemize}
In future studies, we will apply the same reasoning to a given microstructured material, e.g. the Swiss-Cross unit cell \cite{sarhil2024computational}.
Unlike other homogenization methods, we can apply the load very localized within the microstructure.
This requires some modification of the uniform load used here for homogeneous materials, by applying the force at specific points on the central unit cell.
This will give us a direct method to fit the micro-scale elastic tensor $\C_{\rm micro}$ in relaxed micromorphic continua \cite{agn_neff2019static,sky2022primal}.

Additionally, while the analytical solutions of the corresponding cubic symmetry of both the concentrated couple \eqref{eq:CCsolution} and the center of dilatation \eqref{eq:CDsolution} are seemingly more difficult than their isotropic counterpart, they can also be used to numerically find the best approximating cubic counterpart $\C_{\rm cubic}$ of any anisotropic elasticity tensor $\Caniso$.
%
%
%
%
\small
\subsubsection*{Acknowledgments}
Patrizio Neff is grateful for the helpful discussions with Michael A.\ Slawinski (Memorial University of Newfoundland). P.\ Neff acknowledges support in the framework of the DFG-Priority Programme 2256 ``Variational Methods for Predicting Complex Phenomena in Engineering Structures and Materials'', Neff 902/10-1, Project-No. 440935806.
\footnotesize
\section{References}
\printbibliography[heading=none]
\begin{appendix}
%
%
%
\section{Minimal Euclidean/logarithmic distance on $\GLp=\Rp$}\label{app:geodesics}
We want to minimize the Euclidean/logarithmic distance on $\GLp(1)=\Rp$ to find the optimal shear modulus $\mu$ for a given functional form. We start with the logarithmic distance between two positive numbers which is given by \cite{agn_martin2014minimal}
\begin{equation}
	\dist_{\rm geod}^2(x,y)=\left|\log\frac x y\right|^2.\label{eq:distGeodesic}
\end{equation}
We will use this distance measure to best fit a given scalar-valued function $\overline f(r)$ against the function $f_{\rm\mu}(r)$. The goal is to determine the parameter $\mu$.

Hereby, we consider the analytic solution $f_{\rm\mu}(r)=\frac{1}{4\pi\.\mu\.r}$ discussed in Section \ref{sec:solutionConcCouple} with the associated numerical values $(r_i,f_i)$ which are used to give the corresponding interpolation function $\overline f\col(0,r_{\rm max})\to\Rp$. We want to find its closest approximation $f_{\rm\mu}(r)$. For this, we will use the interval $(r_a,r_b)\subset(0,r_{\rm max})$ with $r_a$ as the radius of the inner circle of load and $r_b>r_a$ arbitrary. We consider
\begin{equation}
	F(\mu)=\frac12\.\int_{r_a}^{r_b}\dist_{\rm geod}^2\bigl(\fbar(r),f_\mu(r)\bigr)\,\dr=\frac12\.\int_{r_a}^{r_b}\left[\log\frac{\overline f(r)}{f_\mu(r)}\right]^2\dr\to\min\,.
\end{equation}
We compute
\begin{align}
	F(\mu)=\frac12\.\int_{r_a}^{r_b}\left[\log\frac{\overline f(r)}{f_\mu(r)}\right]^2\dr
	&=\frac12\.\int_{r_a}^{r_b}\left[\log\frac{\overline f(r)}{\frac{1}{4\pi\.\mu\.r}}\right]^2\dr=\frac12\.\int_{r_a}^{r_b}\left[\log(\overline f(r)\.r)+\log(4\pi\.\mu)\right]^2\dr\,.
\end{align}
We search for the minimum by calculating $F'(\mu)=0$, i.e.
\begin{align}
	&&F'(\mu)&=\dd{\mu}\left(\frac12\.\int_{r_a}^{r_b}\left[\log(\overline f(r)\.r)+\log(4\pi\.\mu)\right]^2\dr\right)=\int_{r_a}^{r_b}\left[\log(\overline f(r)\.r)+\log(4\pi\.\mu)\right]\frac{4\pi}{4\pi\.\mu}\dr\overset{!}{=}0\notag\\
	&\iff&0&=\int_{r_a}^{r_b}\log(\overline f(r)\.r)\.\dr+\int_{r_a}^{r_b}\log(4\pi\.\mu)\.\dr\qquad\iff\qquad-(r_b-r_a)\.\log(4\pi\.\mu)=\int_{r_a}^{r_b}\log(\overline f(r)\.r)\.\dr\\
	&\iff&\log(4\pi\.\mu)&=-\frac{1}{r_b-r_a}\int_{r_a}^{r_b}\log(\overline f(r)\.r)\.\dr\hspace{1.77cm}\iff\qquad\mu=\frac{1}{4\pi}\.e^{\displaystyle-\frac{1}{r_b-r_a}\int_{r_a}^{r_b}\log(\overline f(r)\.r)\.\dr}.\notag
\end{align}
We repeat the method for the Euclidean distance
\begin{equation}
	F(\mu)=\frac12\.\int_{r_a}^{r_b}\dist_{\rm euclid}^2\bigl(\fbar(r),f_\mu(r)\bigr)\,\dr=\frac12\.\int_{r_a}^{r_b}\left[\overline f(r)-{f_\mu(r)}\right]^2\dr=\frac12\.\int_{r_a}^{r_b}\left[\overline f(r)-\frac{1}{4\pi\.\mu\.r}\right]^2\dr\to\min\,.
\end{equation}
Again, we search for the minimum by calculating $F'(\mu)=0$, i.e.
\begin{align}
	&&F'(\mu)&=\dd{\mu}\left(\frac12\.\int_{r_a}^{r_b}\left[\overline f(r)-\frac{1}{4\pi\.\mu\.r}\right]^2\dr\right)=\int_{r_a}^{r_b}\left[\overline f(r)-\frac{1}{4\pi\.\mu\.r}\right]\frac{1}{4\pi\.\mu^2\.r}\dr\overset{!}{=}0\notag\\
	&\iff&0&=\int_{r_a}^{r_b}\frac{\overline f(r)}{r}\.\dr-\int_{r_a}^{r_b}\frac{1}{4\pi\.\mu\.r^2}\.\dr\qquad\iff\qquad-\frac{1}{4\pi\.\mu}\left(\frac{1}{r_b}-\frac{1}{r_a}\right)=\int_{r_a}^{r_b}\frac{\overline f(r)}{r}\.\dr\\
	&\iff&\frac{r_b-r_a}{4\pi\.r_ar_b\.\mu}&=\int_{r_a}^{r_b}\frac{\overline f(r)}{r}\.\dr\hspace{2.85cm}\iff\qquad\mu=\frac{\displaystyle r_b-r_a}{\displaystyle 4\pi\.r_ar_b\.\int_{r_a}^{r_b}\frac{\overline f(r)}{r}\.\dr}.\notag
\end{align}
After extensive testing, we concluded that minimizing the logarithmic distance \eqref{eq:distGeodesic} on real numbers does not present any advantage compared to numerically minimizing the Euclidean distance (quadratic minimization).
\end{appendix}
\end{document}